\documentclass[default,iicol]{sn-jnl}% Math and Physical Sciences Reference Style
\jyear{2024}%

\theoremstyle{thmstyleone}
\newtheorem{thm}{Theorem}

\theoremstyle{thmstyletwo}

\theoremstyle{thmstylethree}

\begin{document}

\title[Tilings and grain boundaries]{Colouring monohedral tilings: defects and grain boundaries}


\author{Giedrius Alkauskas}

\affil{\orgdiv{Institute of Computer Science}, \orgname{Vilnius University},  \orgaddress{\street{Naugarduko 24},
\city{Vilnius}, \postcode{LT-03225}, \country{Lithuania}}.  Email: \texttt{giedrius.alkauskas@mif.vu.lt}}



\maketitle

\section{Aperiodicity}
Given a periodic monohedral tiling of the plane. Suppose each tile is assigned one of $N\in\mathbb{N}$ colors, so that neighboring tiles are painted differently. This paper describes a novel phenomenon -- an emergence of \emph{a grain boundary} (1D defect in a 2D material) forced by a given \emph{seed} (a collection of colored tiles).\\   
\indent In the theory of tilings, the pivotal and long awaited-for recent achievement is the discovery of a polygon (one piece, or the \emph{einstein}) which tiles the whole plane, but only non-periodically \cite{strauss, strauss2}. A nice introduction to the subject appeared in the \emph{Intelligencer} as a homage to the legacy of John Conway \cite{radin}. Another paper in this journal \cite{soc-tayl} by Socolar and Taylor describes nicely many aspects in the field: what is a tile; what is a tiling; what is considered to be non-periodic; what are possible matching rules. There is no need to list all the cornerstones of this fascinating story here, specialists can do it better. The reader of this paper will benefit greatly (as did the author) in deeper acquaintance with the works cited. A patch is said to \emph{force} a coloring if there is a unique colored tiling that contains it.\\  
\indent In this note we show a forced full coloring, which is something non-periodic, but surprisingly simple and well-organized. The final object may be thought as being formed by parts of different crystals along with \emph{grain boundaries}; these are interfaces between any two. All results fit into a broader framework of $2$-periodic discrete planar graphs and problem of \emph{unique colorability}, the approach most convenient to work with. It is possible to translate our problem into a tiling setting as described by Taylor and Socolar \cite{soc-tayl} (p. 20-21), and also by Ballier, Durand and Jeandel (\cite{ballier}, Section 2). However, to avoid excessive generality and to make the reading more enjoyable, we stick to the coloring setting.\\
\indent Consider the Taylor-Socolar tile \cite{soc-tayl, soc-tayl2}. It is an example of an \emph{einstein}, with a slight caveat of being non-connected. We emphasize one important feature: there exists a certain 3 tile seed (which does not violate matching rules) which forces the full tiling (\cite{soc-tayl}, p. 26-27). This appears to be in contrast with a \emph{decapod defect} in the Penrose tile setting (\cite{socolar}, p. 240). The latter forces a tiling, too, but the interior of a decapod cannot be filled with Penrose tiles without violating the rules. These two alternatives inspire to ask the following: what happens if several violations to the coloring rule (that is, \emph{defects}) are allowed when dealing with a periodic tiling? As we will see, presence of defects leads to an even broader variety of interesting \emph{polycrystals}.\\

\indent\emph{Dedication.} This paper is dedicated to the memory of Audrius Alkauskas (1978-2023). A physicist, an expert in the field of defects in semiconductors \cite{zhang}. A musician and a writer. My dear twin brother.
\section{The setting}
\begin{figure}
    	\begin{center}
    	\includegraphics[scale=0.32]{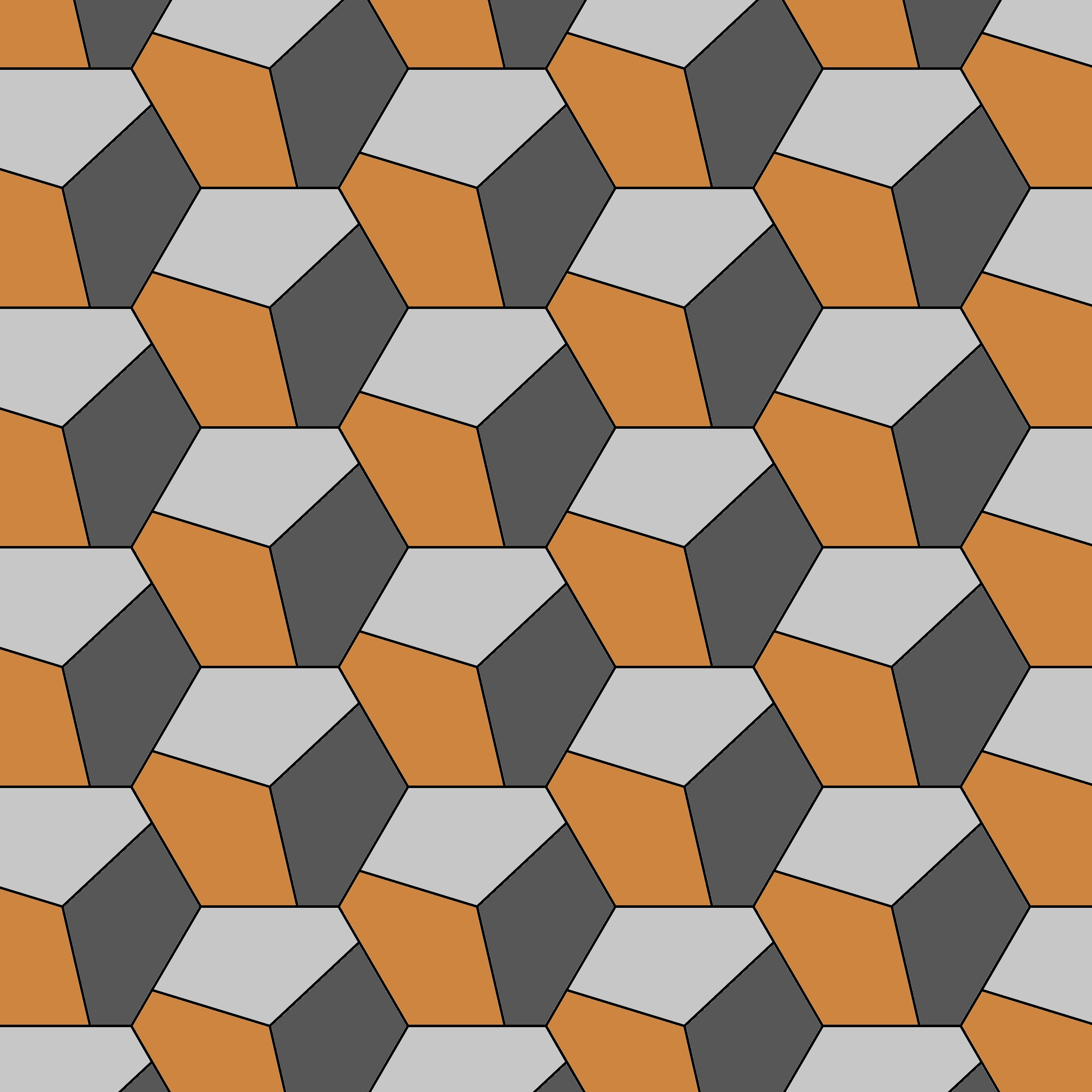} 
    	\caption{Pentagonal tiling of type $3$ (isohedral)} 
    	\label{penta-3}
    	\end{center}
    \end{figure}

\begin{figure}
		\begin{center}
	\includegraphics[scale=0.32]{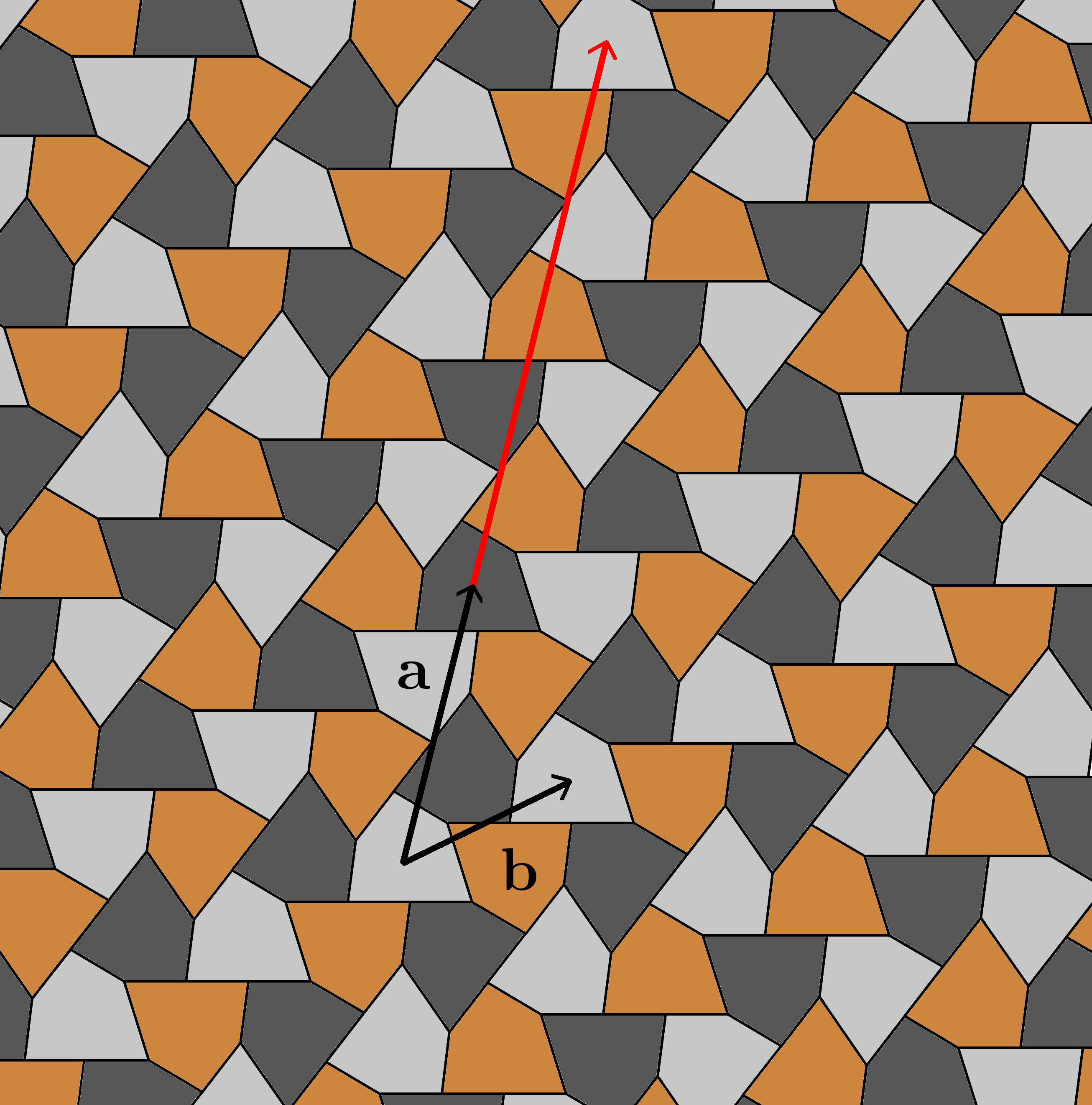} 
	\caption{Pentagonal tiling of type $2$ (isohedral)} 
	\label{penta-2}
		\end{center}
\end{figure}

\begin{figure}
   		\begin{center}
	\includegraphics[scale=0.32]{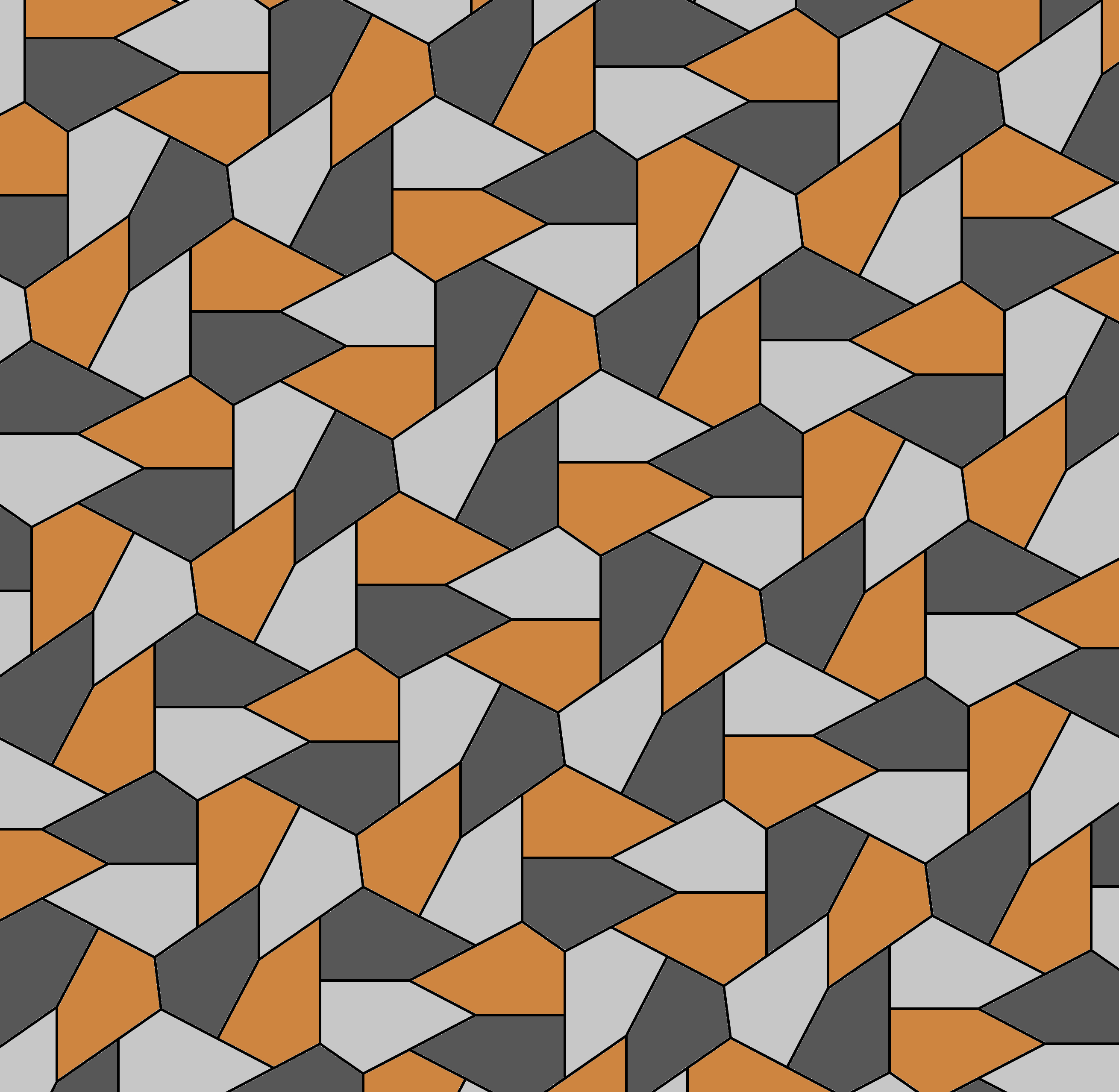} 
	\caption{Pentagonal tiling of type $12$ ($2$-isohedral)} 
	\label{penta-12}
		\end{center}
\end{figure}

Consider (not necessarily edge-to-edge) tessellation of the plane with a single tile (these are called \emph{monohedral}), which at the same time is periodic \cite{shephard}. We say that two tiles are \emph{adjacent} if they share a common point.\footnote{Usually tiles sharing a finite set of points are called \emph{neighbouring}, to distinguish from those sharing a continuous section (or few) of the border (\emph{adjacent}). Yet, the precise definition has only an impact on formation of an adjacency graph. When that is done, the term \emph{adjacent} takes its standard meaning used in graph theory.} Given $N$ colors. We say that a tiling is $N$-colorable, if each tile can be assigned one of $N$ colors in such a way that no two adjacent tiles share the same color.\\
\indent The tiling is said to be \emph{uniquely $N$-colorable}, if such a coloring can be performed (up to permutation of colors) in the unique way. For example, Picture \ref{penta-3} shows pentagonal tiling of type $3$ (discovered by K. Reinhardt) and its unique $3$-coloring. \\
\indent If one wonders whether a unique $N$-coloring of a periodic tiling can be non-periodic, then the answer is, certainly, ``No". Indeed, let us color the whole plane. Let $(\mathbf{a},\mathbf{b})$ be any primitive pair of vectors of periodicity. Shift the tiling with respect $\mathbf{a}$. The new colored tiling, by definition, is obtained from the old one by a permutation of colors. This shows that the tiling is periodic with basis vectors being $k\cdot\mathbf{a}$ and $l\cdot\mathbf{b}$ for certain $k,l\in\mathbb{N}$. Moreover, $k,l$ divide $g(N)$, where $g(N)$ is Landau's function (sequence A000793 in the On-Line Encyclopedia of Integer Sequences), the largest order of an element in the symmetric group $S_{N}$. This sequence  starts as follows:
 \begin{eqnarray*}
 	1, 2, 3, 4, 6, 6, 12, 15, 20, 30, 30, 60, 60, 84, 105.
 	\end{eqnarray*}
 Picture \ref{penta-2} shows another uniquely $3$-colorable pentagonal tiling of type $2$ (also discovered by Reinhardt), where the corresponding vectors are $3\cdot\mathbf{a}$ and $\mathbf{b}$.\footnote{As is clear from Figure \ref{penta-2}, for a primitive pair $(\mathbf{a},\mathbf{a}-\mathbf{b})$, the corresponding vectors would be $(3\mathbf{a},3\mathbf{a}-3\mathbf{b})$. Therefore a pair $(k,l)$ is also not uniquely defined. Finding a maximum value  for $k\cdot l$ (denote this by $t(N)$) over all uniquely $N$-colorable tilings (more generally, $N$-colorable $2$-periodic planar graphs) is a separate interesting problem. We know that $t(N)\leq g^{2}(N)$, which is an equality for $N=3$.} The same phenomenon occurs for type 12 tiling (discovered by M. Rice, Figure \ref{penta-12}). However, the latter is $2$-isohedral (automorphism group of the tiling has $2$ transitivity classes).

\begin{figure}
	\begin{center}
	\includegraphics[scale=0.25]{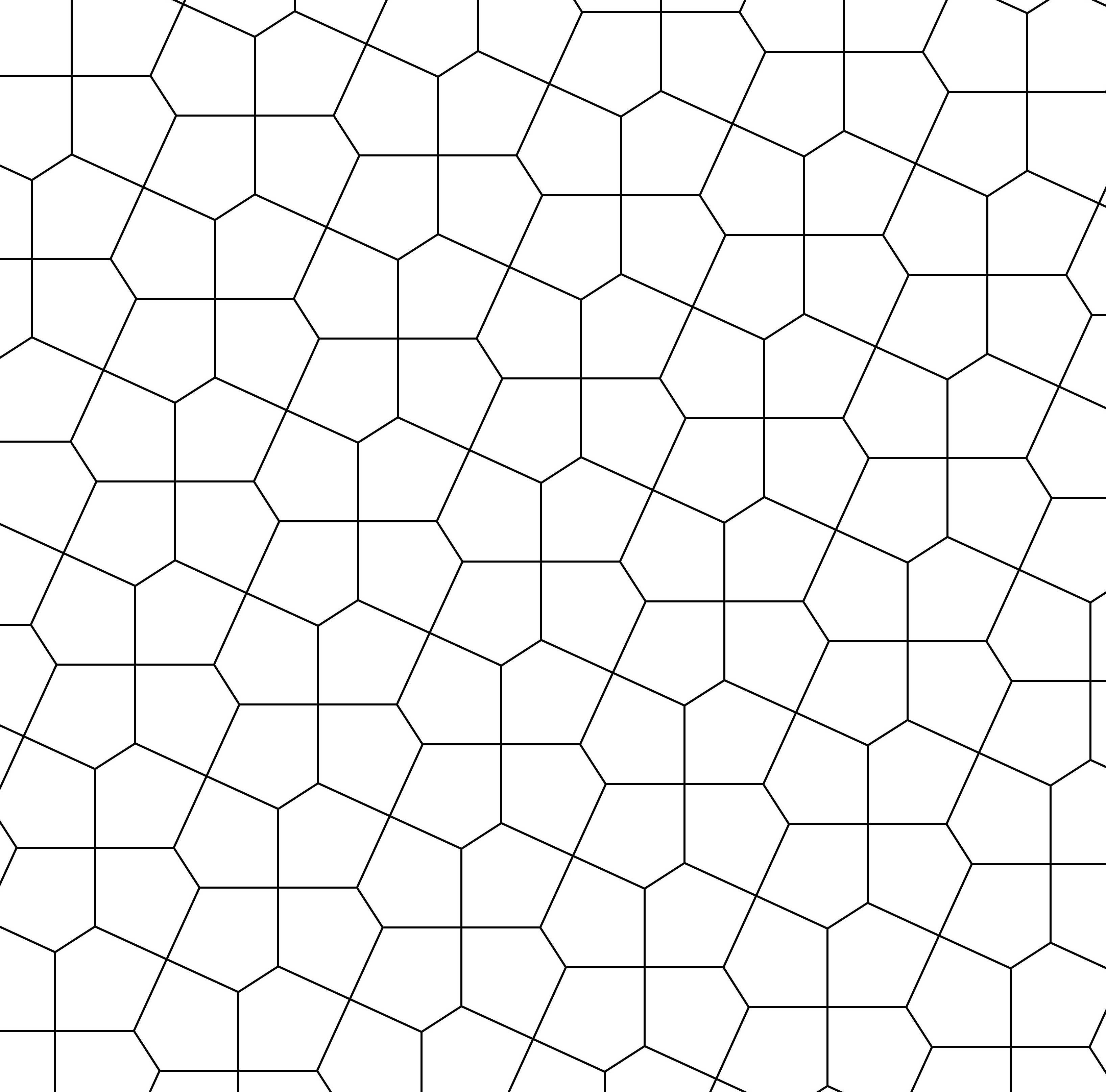} 
	\caption{Pentagonal type $4$ tiling} 
	\label{cairo}
		\end{center}
\end{figure}

 \section{A generic seed}
 \label{generic}
\begin{figure}\begin{center}	
		\includegraphics[scale=0.60]{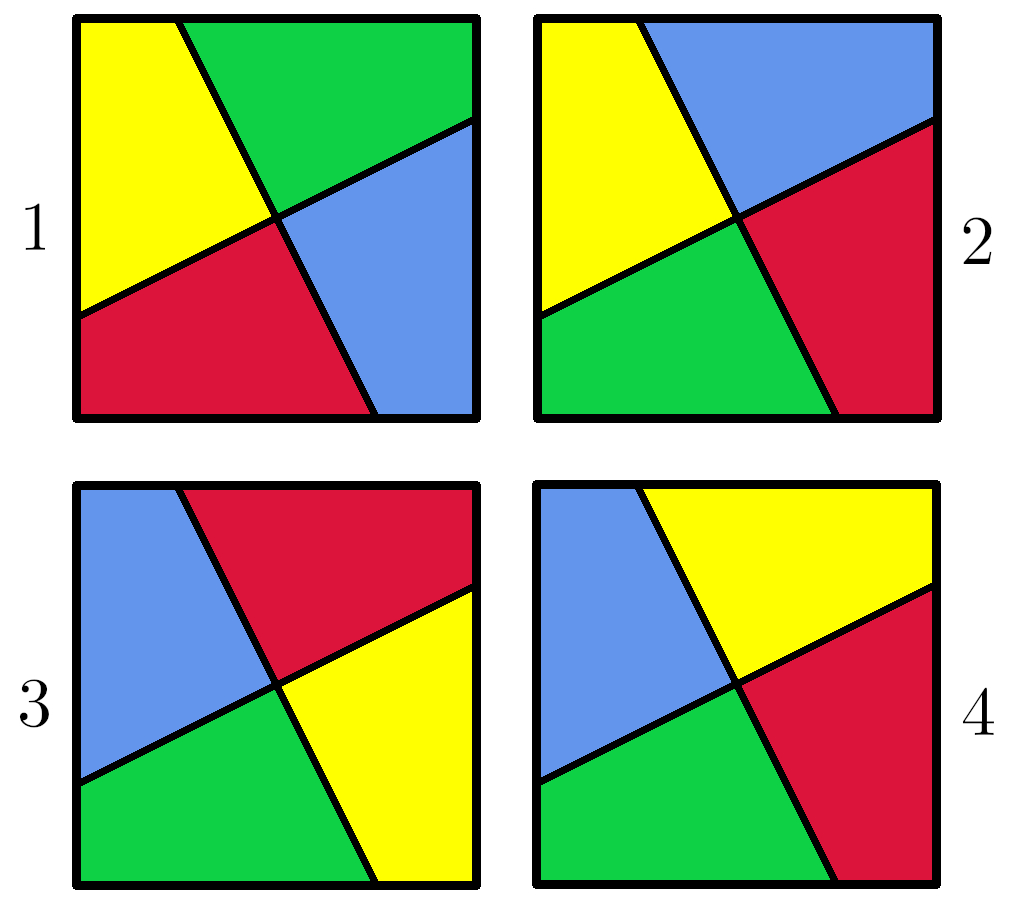} 
		\caption{$4$ (out of $24$) tiling states} 
		\label{unit}
	\end{center}
\end{figure}

Minding the unique colorability property, how one could obtain a uniquely defined non-periodic structure? An analogy with Turing machines and tilings which emulates them (see \cite{strauss}) indicates one possible solution. Indeed, Turing machine is defined by a finite number of states, a finite number of symbols, and a transition function. It starts working when being fed some finite data tape. If we imagine that the unique colorability problem corresponds to an empty (all \texttt{0}'s) tape, then one of the answers to our question (i.e., how to obtain a unique non-periodic structure)  is as follows.\\
\indent Given $N$ colors. Start from a finite number of colored tiles. Call this \emph{a seed}. The only requirement is that this collection does not conflict with the adjacency rule. If the remaining tiles can be colored in exactly $r\in\mathbb{N}$ ways, the initial collection is called \emph{a decent seed} of order $r$. A decent seed of order 1 (which has the unique extension to a full coloring) is called \emph{perfect}. The colored tiling it produces will be also called \emph{perfect}  (we will use terms \emph{full colouring of the tiling} and \emph{full field} interchangeably). \\
\indent Obviously, if one starts from any finite subset of colored tiles in, say, type 2 pentagonal tiling, then either the seed is not contained in any proper coloring, or Figure \ref{penta-2} is obtained with (possibly) colors permuted.\\
\indent  Consider instead a famous Cairo pentagonal tiling (Figure \ref{cairo}). For convenience, we will work with topologically equivalent ``split-chessboard" domino tiling (Figure \ref{1d}), which we label \emph{D-Cairo}. Assume that vertices where $4$ tiles meet form an integer lattice $\mathbb{Z}^{2}$. Fix $N=4$ colors: \textbf{R}ed, \textbf{G}reen, \textbf{B}lue, \textbf{Y}ellow.\\
\indent Let us now divide the whole D-Cairo tiling into anti-clockwise mills (denoted further by A-mills), as shown by black contours in Figure \ref{auga2}. We will heavily employ this division in what is to follow. On the other hand, let as consider $24$ different unit cells, $4$ of which are shown in Figure \ref{unit} (to show every cell, all permutations of 4 colors must be included).
 Simple considerations show that to properly color the D-Cairo board is the same task as to tile the plane with unit cells of $24$ types (``states"), adhering to the following set of restrictions (``patterns"): for any  state, there are $16$ forbidden states to the N(orth), S(outh) (in the picture, 3 cannot be to the South of 1, but 4 can), E(east), and W(est), and 6 forbidden states to each of NE, NW, SE (4 cannot be to the SE of 1), and SW.

\section{A perfect seed}
\label{perr}
\begin{figure}
	\begin{center}
		\includegraphics[scale=0.50]{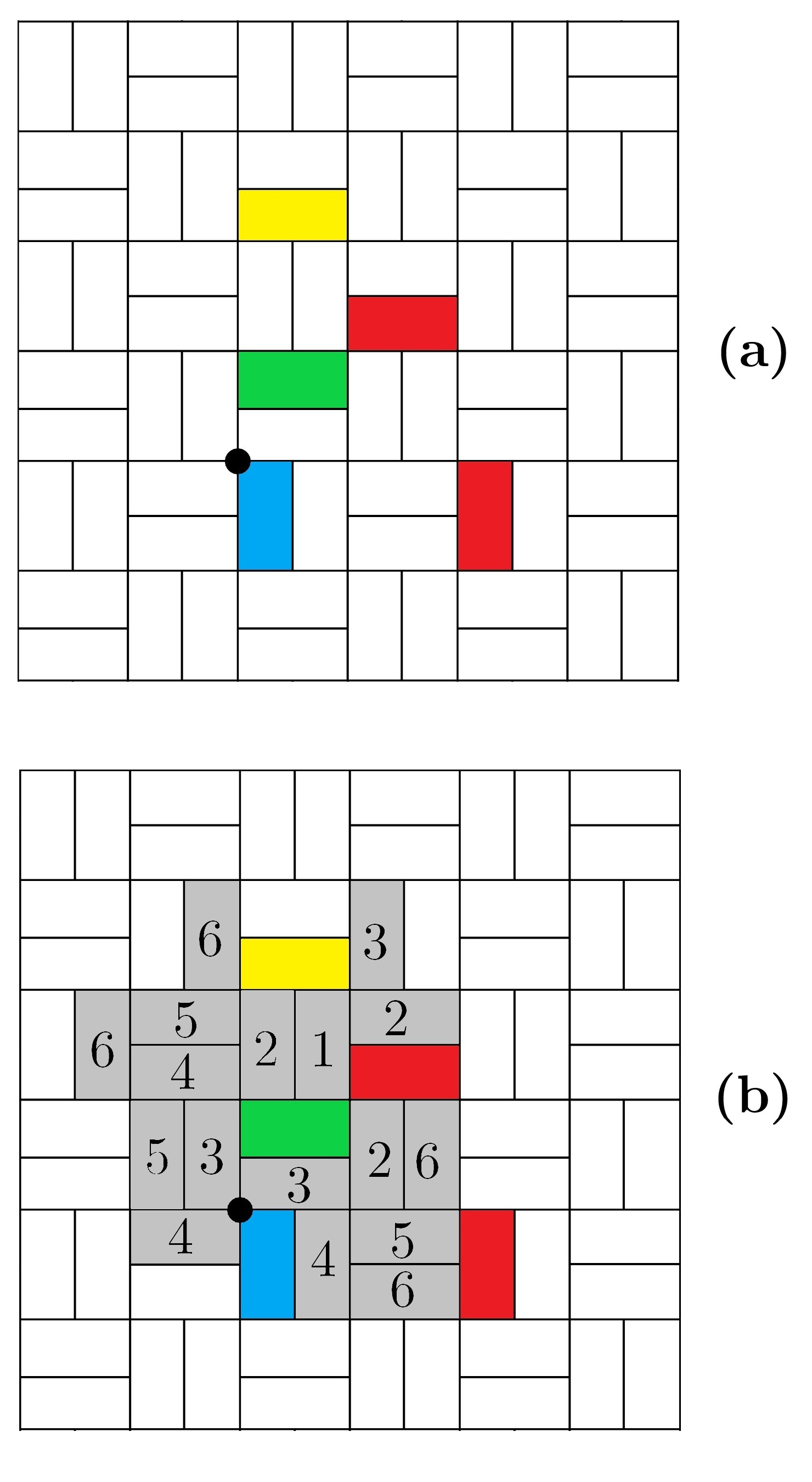} 
		\caption{A perfect seed \textbf{(a)} and its stepwise growth  \textbf{(b)} } 
		\label{auga1}
	\end{center}
\end{figure}
\begin{figure*} 
	\begin{center}
		\includegraphics[scale=0.44]{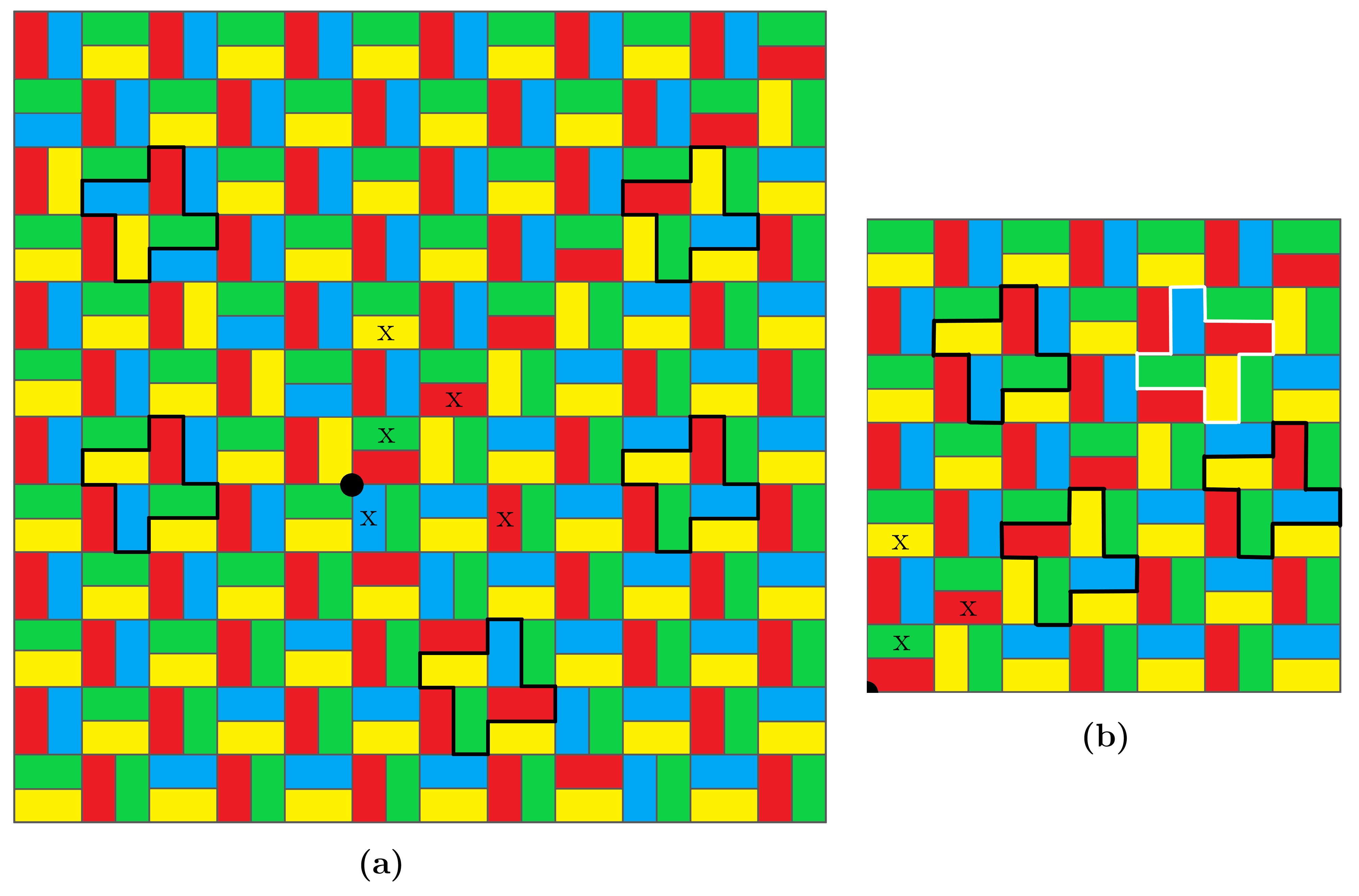} 
		\caption{A central fragment of a perfect field (left) and its upper-right part. The tiles of the seed are labelled by \textrm{x}.} 
		\label{auga2}
	\end{center}
\end{figure*}

\begin{figure}
	\includegraphics[scale=0.43]{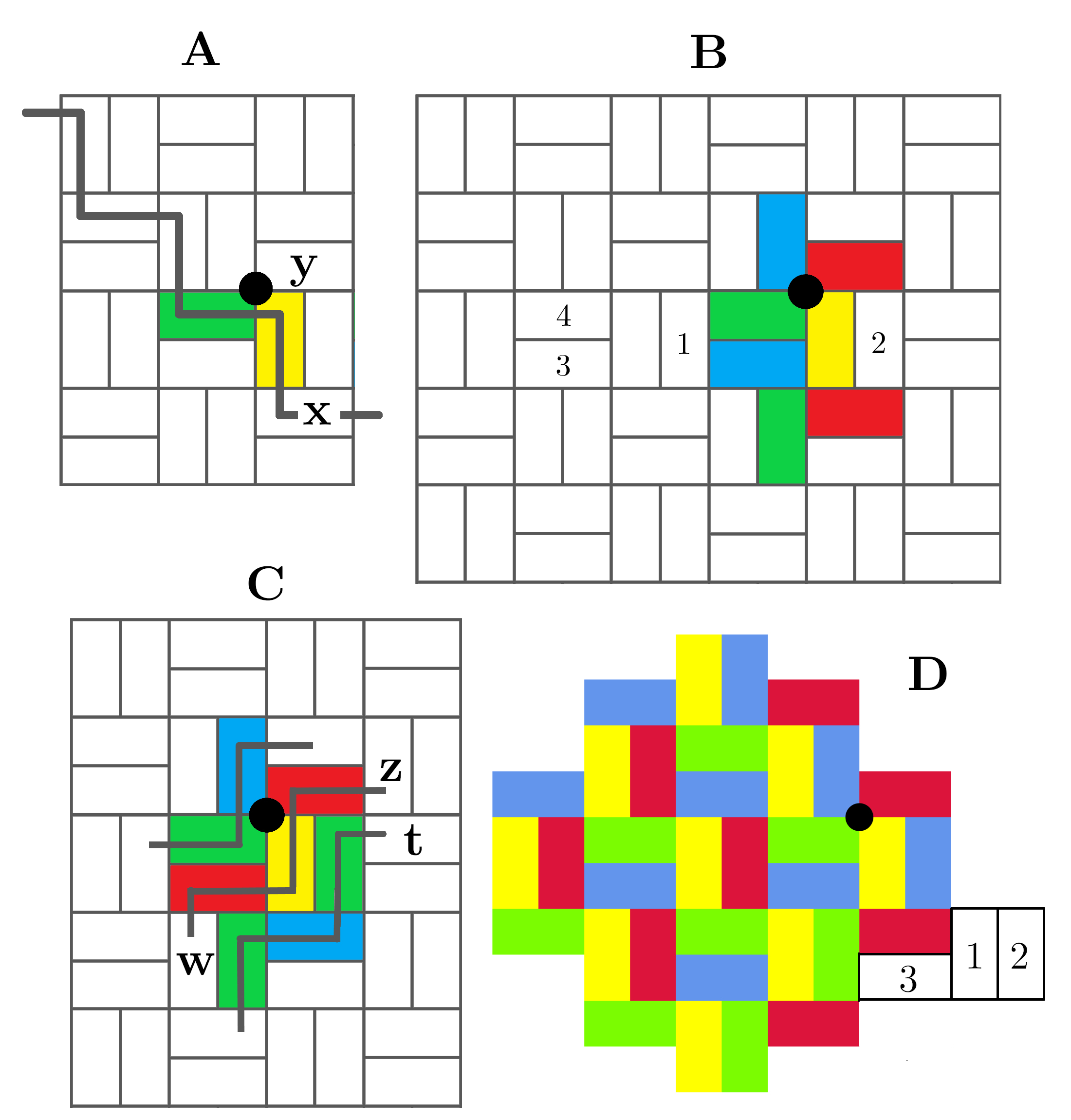} 
	\caption{A configuration in picture \textbf{B} always occurs in a full coloring without type 1 walls. We call this configuration $\mathcal{W}$. It will play a crucial role in proving Theorems \ref{thm1} and \ref{thm2}.} 
	\label{class}
\end{figure}
The fact that perfect seeds exist at first seems surprising.\\
\indent To show this, consider the seed in Fig. \ref{auga1}\textbf{(a)} (a black point -- the center of coordinates -- is for the orientation only). Tile marked as $1$ in \textbf{(b)} should be colored \textbf{B}, since it touches the remaining $3$ colors. The first generation of newly colored tiles thus is a singleton. After that, we can color all the tiles labelled with a number $2$ (the second generation). Then with the number $3$, and so on. Sometimes a new tile to be colored touches 4 tiles of preceding generations, but we never run into a contradiction, and the process continues. The number of tiles (say, $t(n)$) we are able to color at a step $n$ produces the following sequence:
\begin{eqnarray*}
	1, 3, 3, 3, 3, 4, 5, 5, 7, 7, 10, 12, 12, 14, 15, 18, 17,\ldots
\end{eqnarray*}

\begin{figure}
	\includegraphics[scale=0.135]{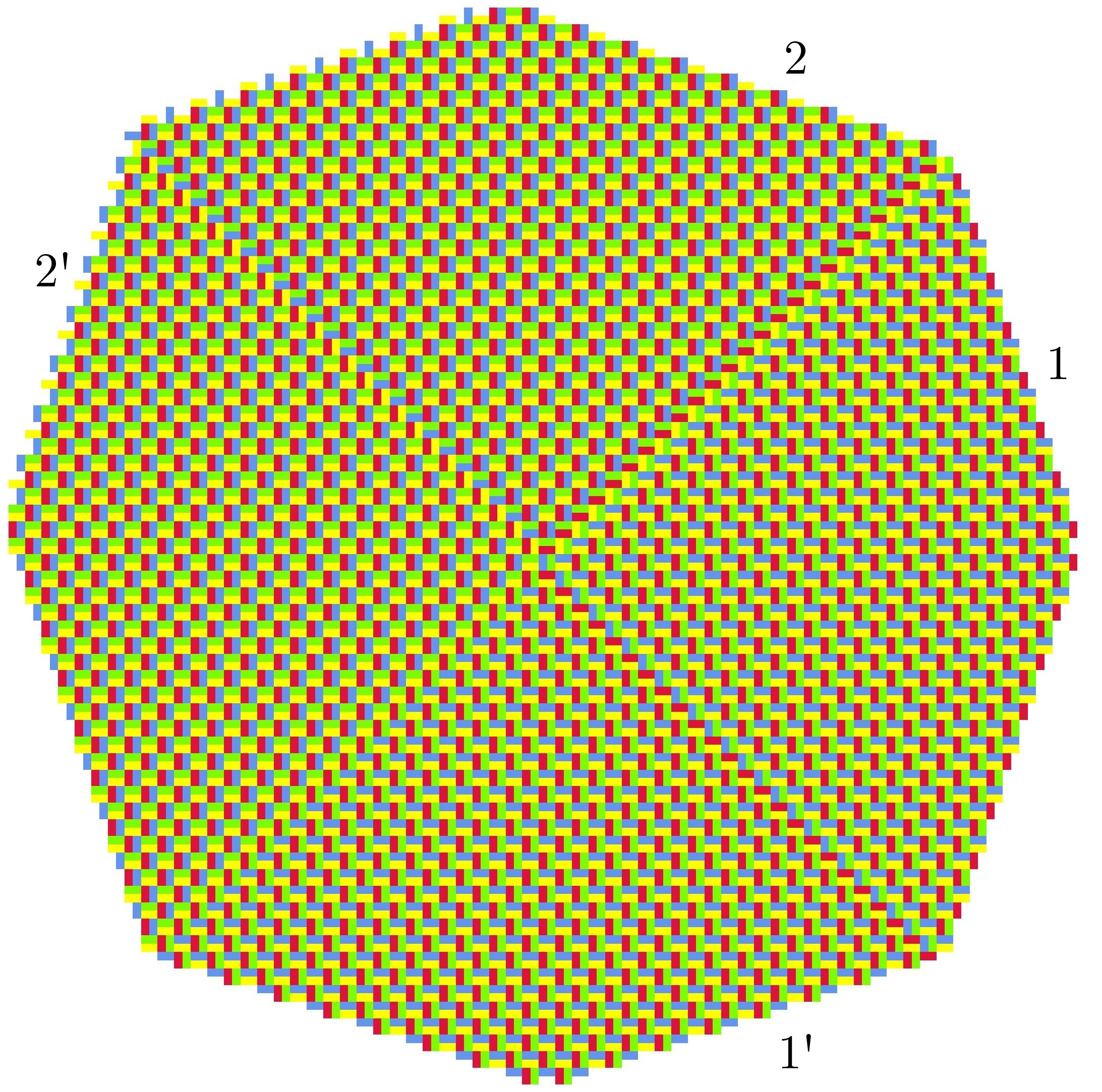} 
	\caption{100 generations produced by the perfect seed in Fig. \ref{auga1}(\textbf{a}).} 
	\label{broken}
\end{figure}

\noindent Its structure becomes apparent when we notice that $t(n+1)-t(n)$ for $n\geq 25$ periodically repeats the length $6$ pattern $( 3, 0, 2, 1, 3, -1)$. Thus for $n\geq 25$, $t(n)=\big{\lfloor}\frac{4n}{3}\big{\rfloor}-3+c(n)$, where $c(n)$ periodically repeats $(-1,1,-1,0,0,1)$.  We also notice that the tiles from a fixed generation form a octagon-like domain, and every tile of the plane will eventually be colored (jumping a bit ahead of ourselves, the first 100 generations are shown in Fig. \ref{broken}). \\
\indent The final result is shown in Fig. \ref{auga2}\textbf{(a)}. The field produced by the seed seems rather periodic, but at certain $1\times 1$ cells the pattern breaks. As already mentioned, let us divide a whole field into A-mills, as shown in Fig. \ref{auga2}\textbf{(b)} by black contours (for a moment, ignore the white clockwise windmill; these will be denoted further by C-mills).\\
\indent It appears that there are only $5$ different types (``states") of mills, shown in Figure \ref{auga2} by black-contours. If we replace each combination with a single color, we arrive to Figure \ref{trys} (top left).  What a surprising picture! Each colored component is a part of a periodic  $2D$ or $1D$ structure (``a crystal"). If we used C-mills instead, we would obtain Figure \ref{trys} (top right), with $6$ distinct colors needed.

This has an appearance of a poly-crystallite material produced by a local assembly rules. The last question we need to answer is the following: does its structure depend on our choice? In other words, if we divide the field into fundamental domains of 4 tiles each in a different fashion (white C-mills as shown in Fig. \ref{auga2}\textbf{(b)}, or into squares $2\times 2$), and define a coloring scheme, will the final result be visually any different? \\
\indent Concerning the coloring of specific domains, sure. Concerning the global picture, no. Indeed, human cognitive system is very good in noticing regularities in visual micro-patters. Essence of this phenomenon is described by psychological \emph{Gestalt laws of grouping}. So, let us color 100 generations of the seed (MAPLE does this job for us). The result is shown in Figure \ref{broken}. Our eyes clearly witness the structure of the material (just zoom the picture and retreat from the screen). In particular, we see that the boundary between $1$ and $2$ is slighly more eminent than the one between $2$ and $2'$. The boundary between $1'$ and $2'$ is faint, while a defect between $1$ and $1'$ is barely noticeable.  
\section{Walls}
\label{wll}
\begin{figure*}
	\includegraphics[scale=0.45]{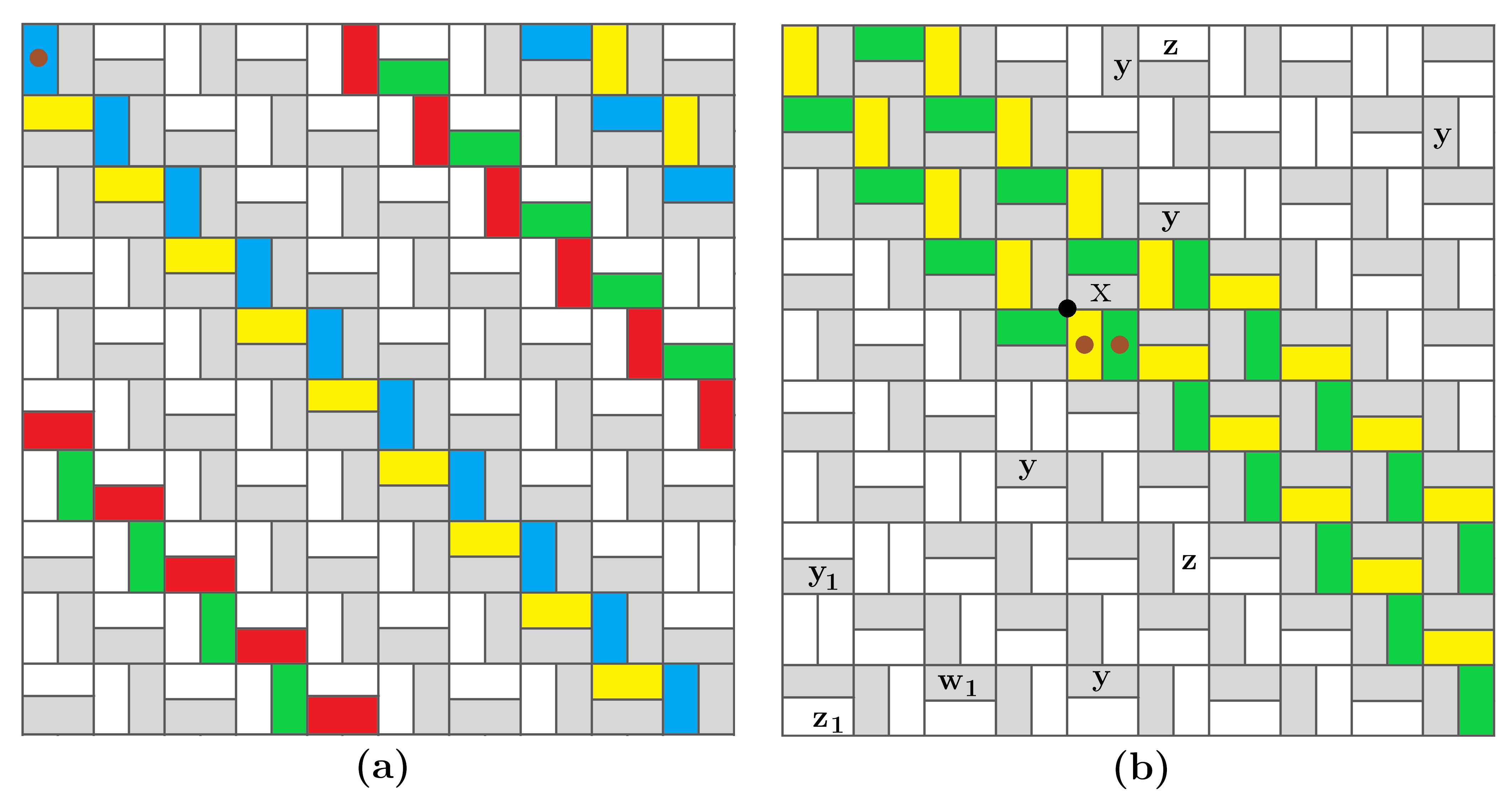}
	\caption{\textbf{(a)} System of type 1 walls in diagonal orientation. Each white staircase should be colored \textbf{B} and \textbf{Y} in alternation, in any of the two ways. The same applies to grey staircases, with colors being \textbf{R} and \textbf{G}.\\
		 \textbf{(b)} System of type 2 walls. Each while staircase (type 2 wall) should be colored with \textbf{Y} and \textbf{G} in the unique way. Each grey staircase (a trail) should be colored \textbf{R} and \textbf{B} in alternation, in any of the two ways.} 
	\label{1d}
\end{figure*}

As will soon become clear, the reason why a particular seed can be extended to a full coloring in infinitely many ways is that it does not forbid a formation of a system of walls. \emph{A wall} is a diagonal chain (a staircase) of tiles, colored alternately in two colors. Walls can be of two basic types. \emph{System of walls} is a very specific full coloring of a tiling. Such a system is, consequently, also of two types. Both types can be in diagonal (see Fig. \ref{1d}) or in  antidiagonal orientation.\\  
\indent In order to show the first prototype, suppose a diagonal chain in two colors \textbf{B}\textbf{Y} extends through the whole tiling, infinitely to both ends, as shown in Fig. \ref{1d}(\textbf{a}) (a staircase with a brown dot). 
\begin{figure*}
	\includegraphics[scale=0.65]{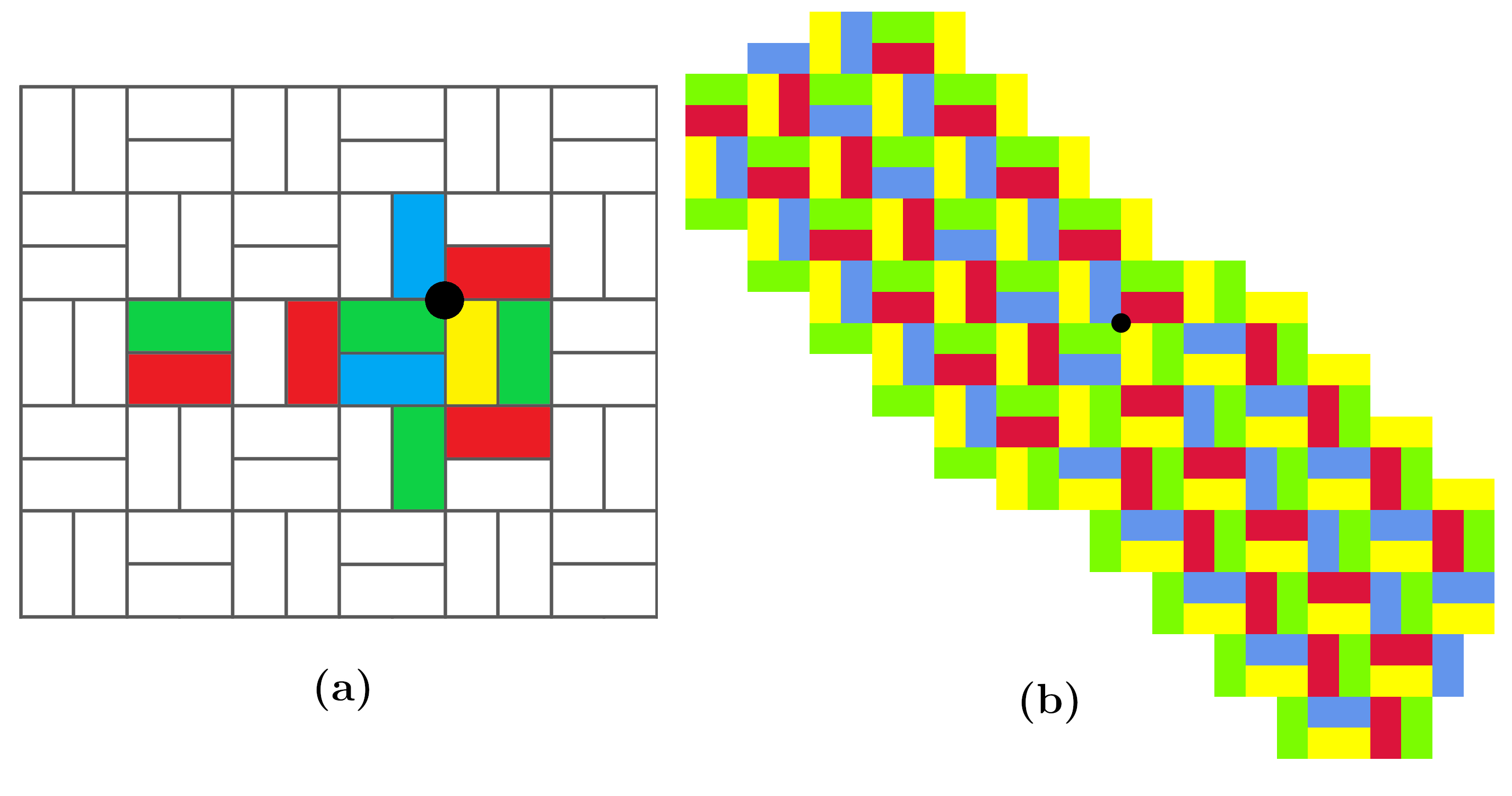}
	\caption{The seed \textbf{(a)} leads \textbf{(b)}, a formation of four \textbf{Y}\textbf{G} type 2 walls and three \textbf{R}\textbf{B} trails} 
	\label{antraip}
\end{figure*}
We call this chain a \emph{type 1 wall}. It is periodic. Another type 1 wall must be built immediately below (and above it) to comply with the adjacency rule. For its coloring we have two choices: \textbf{R} and \textbf{G} in alternation, or the other way round. It is now obvious that each such coloring is in 1-to-1 correspondence with a doubly infinite sequence of \texttt{0}'s and \texttt{1}'s. As is well known, this has the power of continuum.\\
\indent  In Section \ref{all-poly} we will deal with finite connected subsets of type 1 walls. Call this formation a \emph{fence}, if it contains at least $2$ colored tiles. For example, do not extend all four type 1 walls shown in Figure \ref{1d} to infinity but rather consider the picture as it is. Then the picture depicts 4 diagonal fences. Their lengths are, respectively, $9$, $19$, $11$, $5$.  By definition, two tiles shown in Fig. \ref{class}\textbf{A} form a fence of length $2$. More generally, every tile in the full field belongs to one diagonal and one antidiagonal fence. Any fence can be extended until it encounters a tile of the third color. If it can be extended indefinitely to both ends, we have a type 1 wall. \\ 
\indent In order to introduce the second basic prototype, let us consider the finite seed given in Fig. \ref{antraip}\textbf{(a)}. If we color $25$ generations produced by this seed, we obtain the picture \textbf{(b)}. This is a new type of behavior: infinitely many tiles can be uniquely colored (soon we will see that many more have a pre-defined color!), and they are contained between two identical formations, colored alternatively \textbf{Y}\textbf{G}. Each of them extend infinitely to both directions. However, these are NOT type 1 walls, but something slightly different. One of these formations is shown the Fig. \ref{1d}\textbf{(b)}; namely, the staircase containing two dotted tiles. Call this \emph{type 2 wall}. It is non-periodic. This is the main difference with type 1, since type 2 wall has two distinguished colored tiles, exactly those with brown dots. \\
\indent  Let us take a closer look at type 2 wall in Fig. \ref{1d}(\textbf{b}). The tile \textbf{x} can be either \textbf{R} or \textbf{B}. In each case this forces a coloring of a diagonal strip of tiles immediately above (shown by light grey) which we label \emph{a trail}. After this one is colored, we witness the formation of another type $2$ wall of exactly the same pattern as the initial one, independently of our choice for the coloring of \textbf{x}. In other words, the second type 2 wall immediately above is predestined, while a trail between them can be colored in two ways. The process continues infinitely to both sides. In this way the whole field is divided into type 2 walls and trails. Thus, tiles labelled by \textbf{y} can each be colored either \textbf{R} or \textbf{B}. This forces the coloring of the full trail to which each of these tiles belong. On the other hand, tiles labelled by \textbf{z} must be colored in the unique way. It is now clear that the seed in Fig. \ref{antraip}\textbf{(a)} can be completed to a full field in a continuum of ways. \\
\indent Moreover, any patch which contains it as a subset either has no extensions at all, or can also be extended in a continuum of ways. For example, suppose these additional tiles are $\mathbf{y}_{1}$, $\mathbf{z}_{1}$, $\mathbf{w}_{1}$; see Fig. \ref{1d}(\textbf{b}). If $\mathbf{z}_{1}$ is \textbf{G}, \textbf{R} or \textbf{B}, there is no extension. This holds because $\mathbf{z}_{1}$ belongs to a type $2$ wall and must be colored in the unique way (in this case, \textbf{Y}) in order an extension to exist. Similarly, $\mathbf{y}_{1}$ and $\mathbf{w}_{1}$ belong to the same trail. Hence, unless $\mathbf{y}_{1},\mathbf{w}_{1}$ are colored \textbf{R}, \textbf{B} or  \textbf{B}, \textbf{R},\footnote{Assign $0$ to any tile in the trail, and afterwards assign a corresponding integer to its every tile. What is essential of tiles belonging to this extended seed is their parity.} there is no extension. If, however, say, $\mathbf{y}_{1},\mathbf{w}_{1}$ have colors \textbf{R}, \textbf{B}, this forces a unique coloring of the trail, and there are a continuum of extensions.\\
\indent It is now clear that, given a finite patch, the first thing one needs to check is whether it belongs to a system of type 1 or type 2 walls, either diagonal or antidiagonal. If this is the case, there exists a continuum of extensions. In the next two Sections we will see that otherwise the number of extensions is always finite (possibly, none at all).

\section{All polycrystals} 
\label{all-poly}
\begin{figure*}
	\begin{center}
		\includegraphics[scale=0.70]{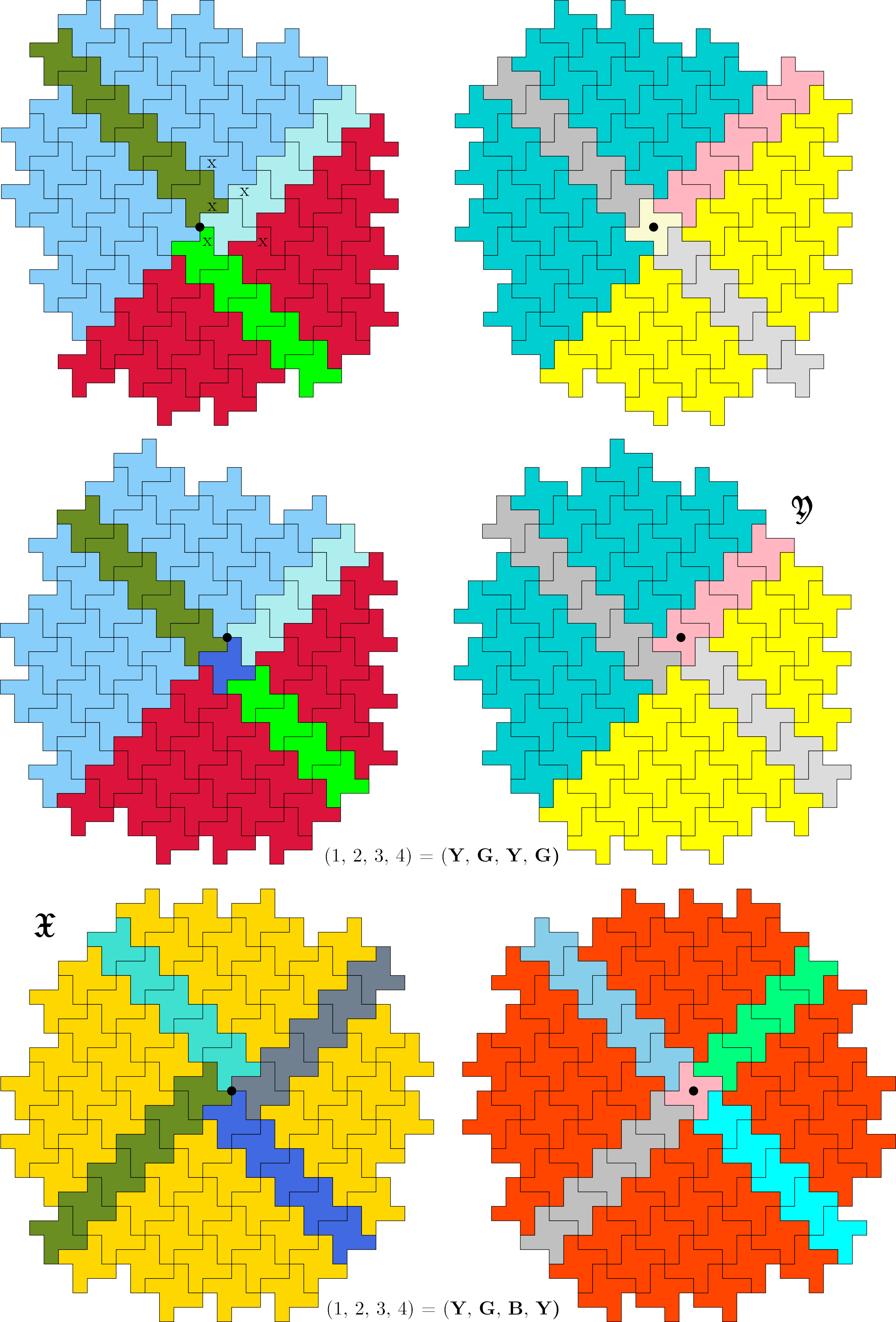} 
		\caption{Fields of perfect seeds divided into windmills} 
		\label{trys}
	\end{center}
\end{figure*} 
Many natural question at this point arise. For example: are there any other essentially distinct examples than the one in Fig. \ref{broken}? For the starters, note that the symmetry group of a D-Cairo tiling (denoted in crystallography by $p4g$ or by 4*2 in the orbifold notation) includes the reflection with respect to the line $x=\frac{1}{2}$. However, when we collect all $4$ dominoes around every vertex $(x,y)\in\mathbb{Z}^{2}$, $x+y\text{ (mod }2)=1$ (``odd vertices", as opposed to ``even vertices") into A-mills, the reflection symmetry is lost: such transformation caries an A-mill into a C-mill. Yet, this particular grouping is our subjective choice, it is not present in a tiling itself. This remark has the following consequence: when we classify all perfect fields, both versions of mills must be considered. It is natural to call colored tilings \emph{equivalent}, if they coincide after a (possible) symmetry of the tiling and permutation of colors. With that said, we can state our first result. 
\begin{thm}Suppose a finite $4$-color seed in the D-Cairo tiling forces the full coloring. Then the full field is equivalent to one of two colored tilings: $\mathfrak{X}$ or $\mathfrak{Y}$. When divided into mills, these have the appearance shown in Fig. \ref{trys}.
	\label{thm1}
	\end{thm}
To prove this, suppose we have a perfect seed. Let us color the whole field. At every vertex (even or odd) where $4$ tiles meet, all $4$ colors must be present. We will choose one such vertex as a starting point having the following additional criteria in mind. \\
\indent As we have seen in Section \ref{wll}, our field cannot contain a type 1 wall (it also cannot contain a type 2 wall; this fact will be used later). Choose any even vertex. Let bottom-left and bottom-right colors at it be \textbf{G} and, respectively, \textbf{Y} (Figure \ref{class}\textbf{A}). These two tiles form a diagonal fence. We go along this fence (shown by grey staircase) in the direction where it terminates. This must happen to one side or another, otherwise a fence belonged to a type 1 wall. Termination occurs when a third color in encountered. Say, at a tile \textbf{x} (recall that, if needed, we can perform a permutation of colors, a shift of a tiling, and a reflection with respect to a line $x=\frac{1}{2}$).\\
\indent If the colors at \textbf{x} and \textbf{y} are the same, we obtain the seed shown in Fig. \ref{class}\textbf{B}. Call this and equivalent configurations as $\mathcal{W}$. If the colors are distinct, the obtained seed is shown in \textbf{C}. However, in the latter case we witness a formation of three parallel fences, this time antidiagonally oriented (shown by grey staircases). Let us go along the middle \textbf{R}\textbf{Y} fence in the direction where it terminates. If this is a tile \textbf{z}, then it must be colored \textbf{B}. A tile \textbf{t} is now \textbf{Y}, and (after a rotation and permutation of colors) we recognize $\mathcal{W}$ as a subset. If a fence terminates at \textbf{w}, then this tile is colored \textbf{B}. Five more tiles then have a forced color, and we again witness $\mathcal{W}$ as a subset. If the fence in question terminates at some more distant tile, it drags two surrounding fences alongside (their coloring is forced), and we are back to the setting just discussed. Therefore, without loss of generality, we can take the configuration $\mathcal{W}$ as a starting point -- it is always to be found in a field which is not a system of type 1 walls. Tile \textbf{Y} in it is exactly the one where two distinct diagonal, as well as two distinct antidiagonal fences meet.  \\
\indent Choose $4$ tiles (marked $1,2,3,4$ in the figure) so that there is no conflict with the adjacency rule. For example, we can take colors $(\mathbf{Y},\mathbf{G},\mathbf{B},\mathbf{G})$. In total, there are $2\times 2\times 4\times 3=48$ possibilities.  MAPLE code works on such input \cite{maple}. If eventually a cell neighbouring $4$ distinct colors is encountered, the program halts with a FAULT. If the program can run on infinitely, it produces a pre-given number of generations.\\
\indent  In fact, it is enough to check only $24$ cases. To see this, note that $\mathcal{W}$ is invariant under the following transformation $\mathcal{T}$: first, reflect the board with respect to the line $y=-\frac{1}{2}$, second, swap \textbf{B} and \textbf{G}. Since tile 2 is necessarily \textbf{B} or \textbf{G}, without loss of generality, we may consider 2 painted in \textbf{B}.  \\
\indent Here is the outcome of these computations. The program halts in $18$ out of $24$ cases. In $4$ other cases we witness the growth of a polycrystal very similar to the one described in Section \ref{perr}. However, a closer look reveals some differences. In order to understand the situation better, an extension of the program was added. The latter automatically subdivides the board into A-mills (and, independently, C-mills) and re-colors it. In general, $24$ colors are needed. However, if we accept the metaphor of a ``mill", we may consider cyclic permutations of colors as being close variants of one another. This suggests the following scheme. Six color groups, each consisting of four hues, are chosen: yellow-gold, green-olive, grey-silver, red-pink, cyan-turquoise,  blue-violet. Permutations of $4$ colors belong to the same group if they are cyclic. This convention applies to A-mills. For C-mills, we imagine it reflected with respect to the vertical line. Now it is A-mill and it has its color.\\
\indent The rest is done by the program. The choice $(1,2,3,4)=(\mathbf{Y}, \mathbf{G}, \mathbf{Y}, \mathbf{G})$\footnote{Here and in the next example tile 2 in Fig. \ref{class}\textbf{B} is painted \textbf{G}. Under the transformation $\mathcal{T}$ this corresponds to a quadruple $(\mathbf{Y}, \mathbf{B}, \mathbf{B}, \mathbf{Y})$.} leads to the picture shown in Figure \ref{trys} (middle). Analogous picture (in case of necessity we swap A-mill with C-mill) occurs in two more cases. We label the C-mill configuration as $\mathfrak{Y}$. On the other hand, the choice  $(1,2,3,4)=(\mathbf{Y}, \mathbf{G}, \mathbf{B}, \mathbf{Y})$ leads to the picture shown at the bottom. We label the A-mill version as $\mathfrak{X}$. Its pleasingly symmetric, and the whole polycrystal consists of one single crystal with four $1D$ defects. Its C-mill version shows that a central ``molecule" (light pink) is a rotation of the main (orange red) molecule. \\
\indent Next, out of $24$ cases, $1$ leads to the formation of type $2$ wall and so must be discarded. If we chose tile 2 to be \textbf{G} rather than \textbf{B}, this distinguished seed would be exactly the one shown in Fig. \ref{antraip}(a). \\
\indent We are left to consider the final choice $(1,2,3,4)=(\mathbf{R}, \mathbf{B},\mathbf{B}, \mathbf{G})$. In this case the program does not halt, but it colors only the tiles shown in Figure \ref{class}\textbf{D}. Now, all 9 possible choices to color tiles 1, 2, 3 lead either to a type 2 wall, to a contradiction, or a polycrystal already encountered.\\
\indent One final remark. If an obtained A-mill or C-mill picture has patterns that coincide with those of $\mathfrak{X}$ or $\mathfrak{Y}$ after (a possible) shift, rotation and reflection, it does not automatically imply that full fields are equivalent, unless one goes deeper into the coloring scheme and interrelations among all $24$ colors. Due to this caveat, a program and a seed list \cite{maple} contain few lines which allow to permute our $4$ basic colors and rotate a field. This helps in demonstrating that a polycrystal which, say, looks like $\mathfrak{Y}$, is in fact equivalent to $\mathfrak{Y}$. This completes the proof.
\section{Decent seeds and defects}
\label{decent}
\begin{figure}
	\begin{center}
		\includegraphics[scale=0.22]{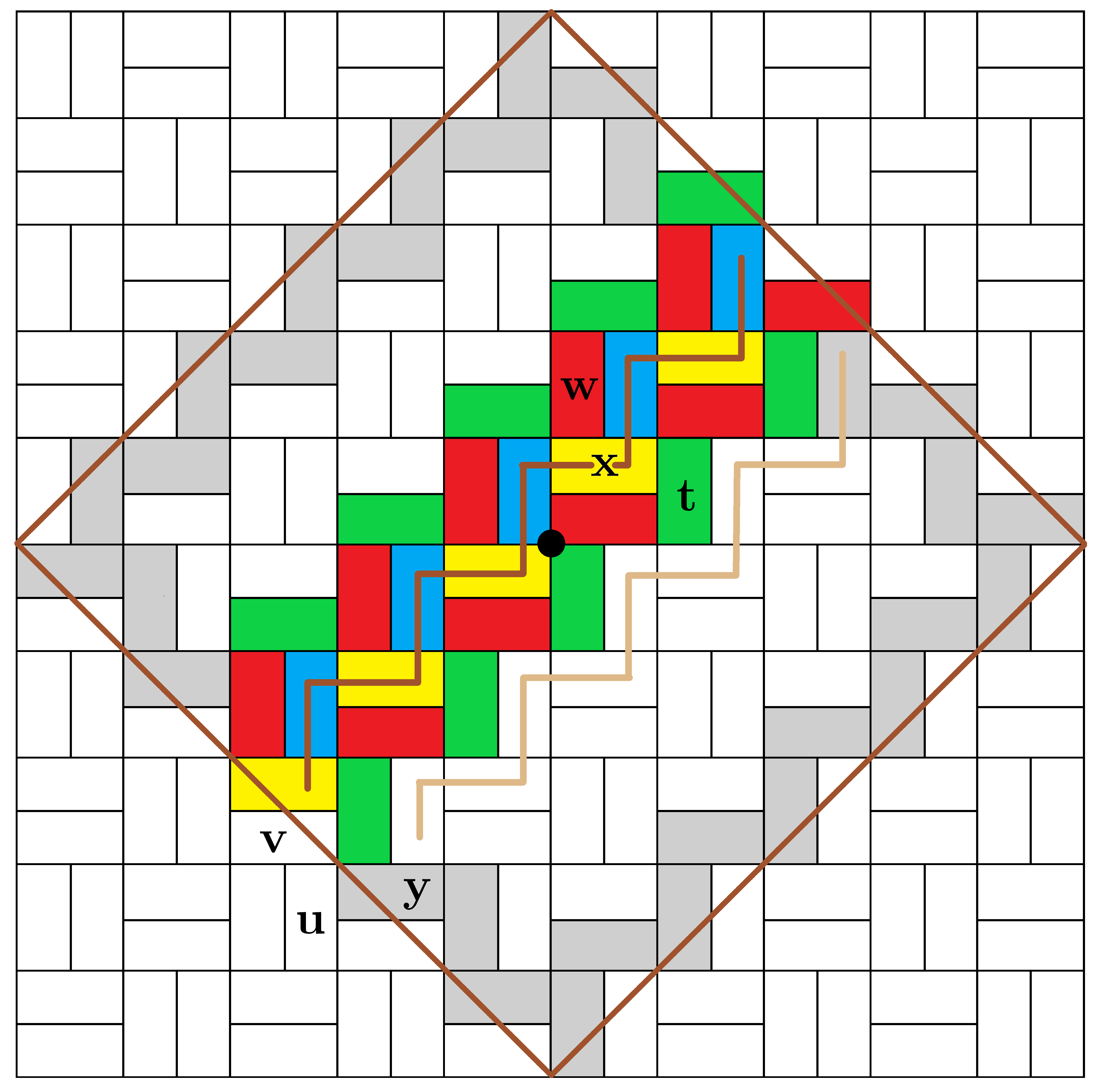} 
		\caption{A supplementary construction: a set $\mathcal{S}_{5}$ consists of grey (boundary) tiles and those fully inside a brown square.}
		\label{saib} 
	\end{center}
\end{figure}

\begin{figure}
	\begin{center}
		\includegraphics[scale=0.43]{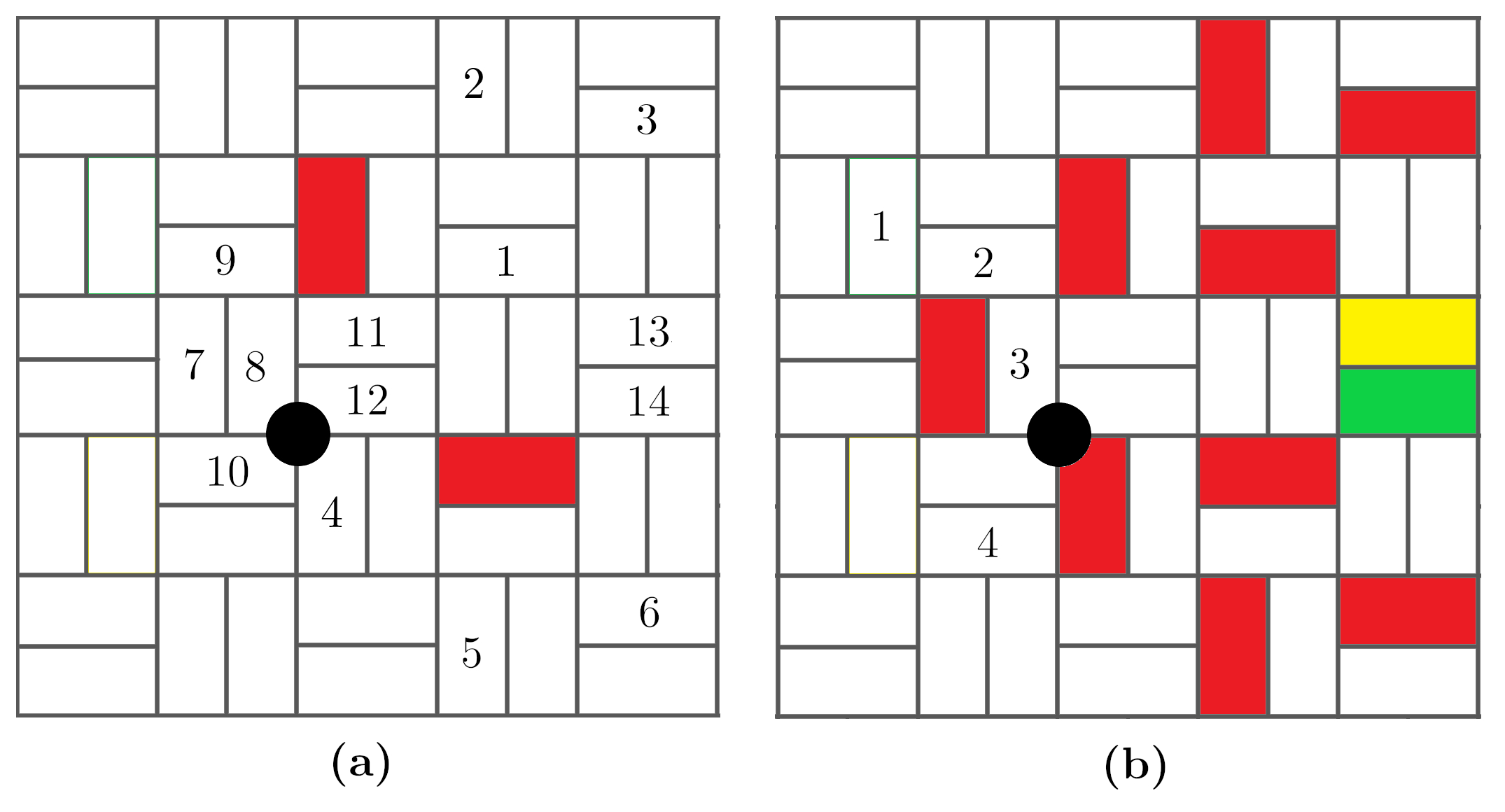} 
		\caption{A decent seed (a), and one of 6 its possible extensions (b)} 
		\label{almost}
	\end{center}
\end{figure}

Having dealt with perfect fields, we have all the tools needed to prove the second main result. 
\begin{thm}
The cardinality of extensions of any given seed is either finite or a continuum.\\		
\indent The 2-tile seed shown in {\rm Fig. \ref{almost}\textbf{(a)}} is decent of order $36$. There exist decent seeds of orders $1,2,3,4,6,8,12$, and also of order exceeding any given bound.
\label{thm2}
\end{thm}
Let $M\in\mathbb{N}$. Consider all tiles which either fully, or the the bigger part of them lie inside the square $\vert x\vert+\vert y\vert\leq M$. Call this set $\mathcal{S}_{M}$ (Fig. \ref{saib} depicts $\mathcal{S}_{5}$).   

\begin{figure}
	\includegraphics[scale=0.60]{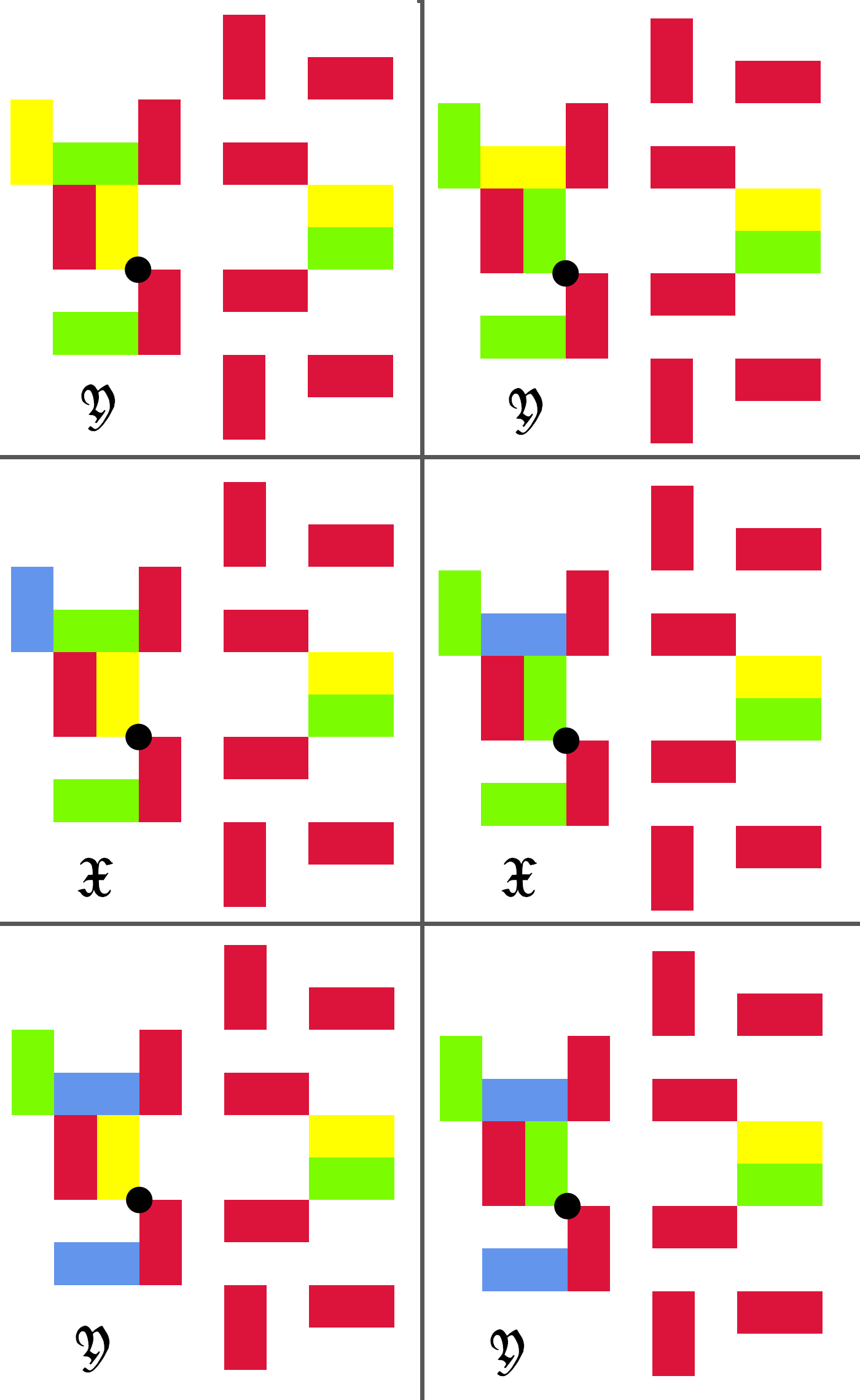} 
	\caption{Perfect seeds, possible extensions a decent seed in Fig. \ref{almost}(\textbf{b}).
		Permutations of \textbf{G}, \textbf{Y}, \textbf{B} give $36$ extensions of the seed in Fig.\ref{almost} (\textbf{a}).} 
	\label{almost2}
\end{figure}

\begin{figure*}
 	\begin{center}
 		\includegraphics[scale=0.52]{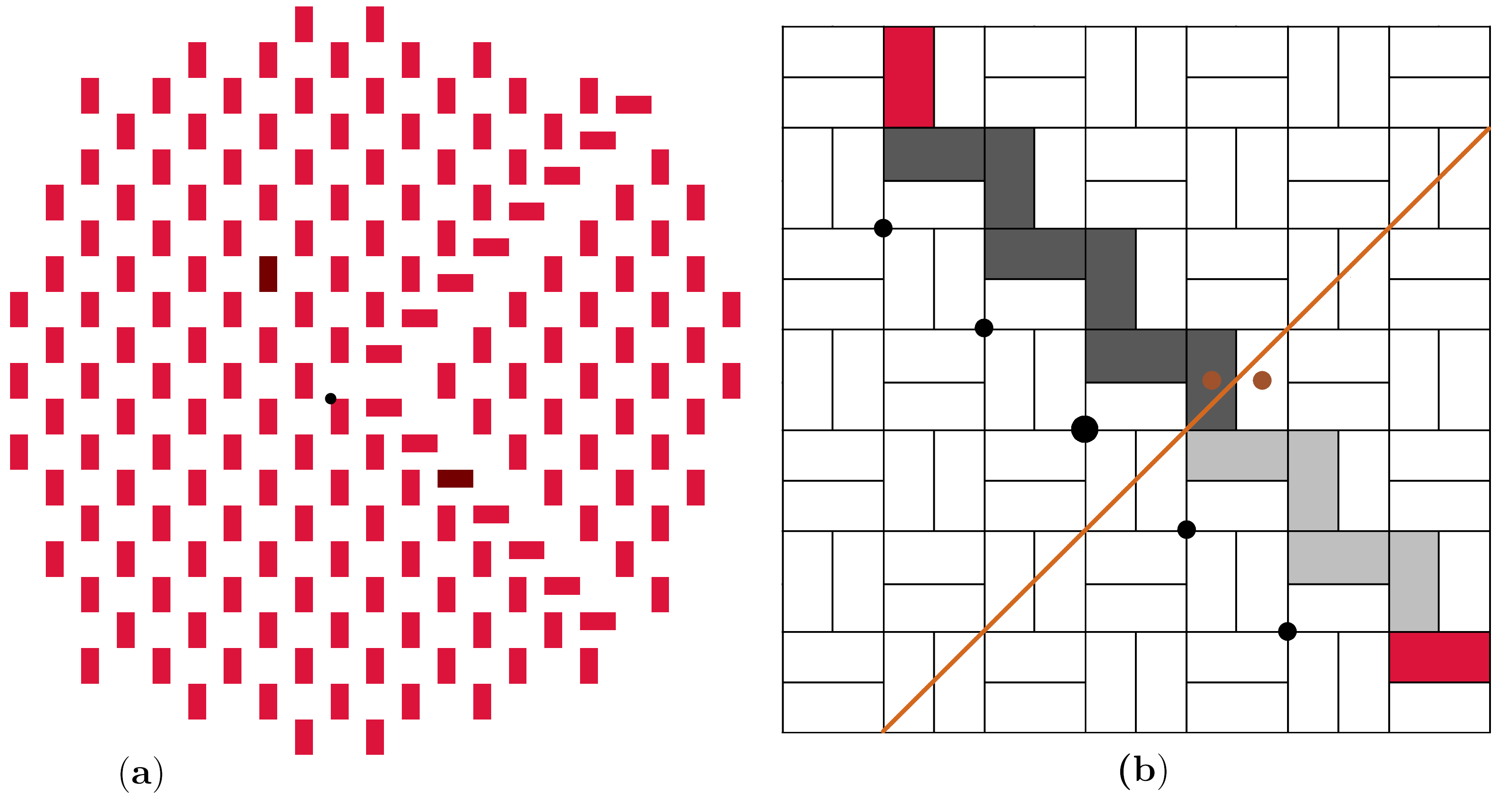} 
 		\caption{(\textbf{a}) \textbf{R} tiles of all $36$ polycrystals as described by Fig. \ref{almost2}; two red tiles are marked with dark-red. These two tiles are shown in (\textbf{b}); it is a decent seed or order of at least $180$.}
 		\label{raud} 
 	\end{center}
 \end{figure*}

\begin{figure*}
 	\includegraphics[scale=0.72]{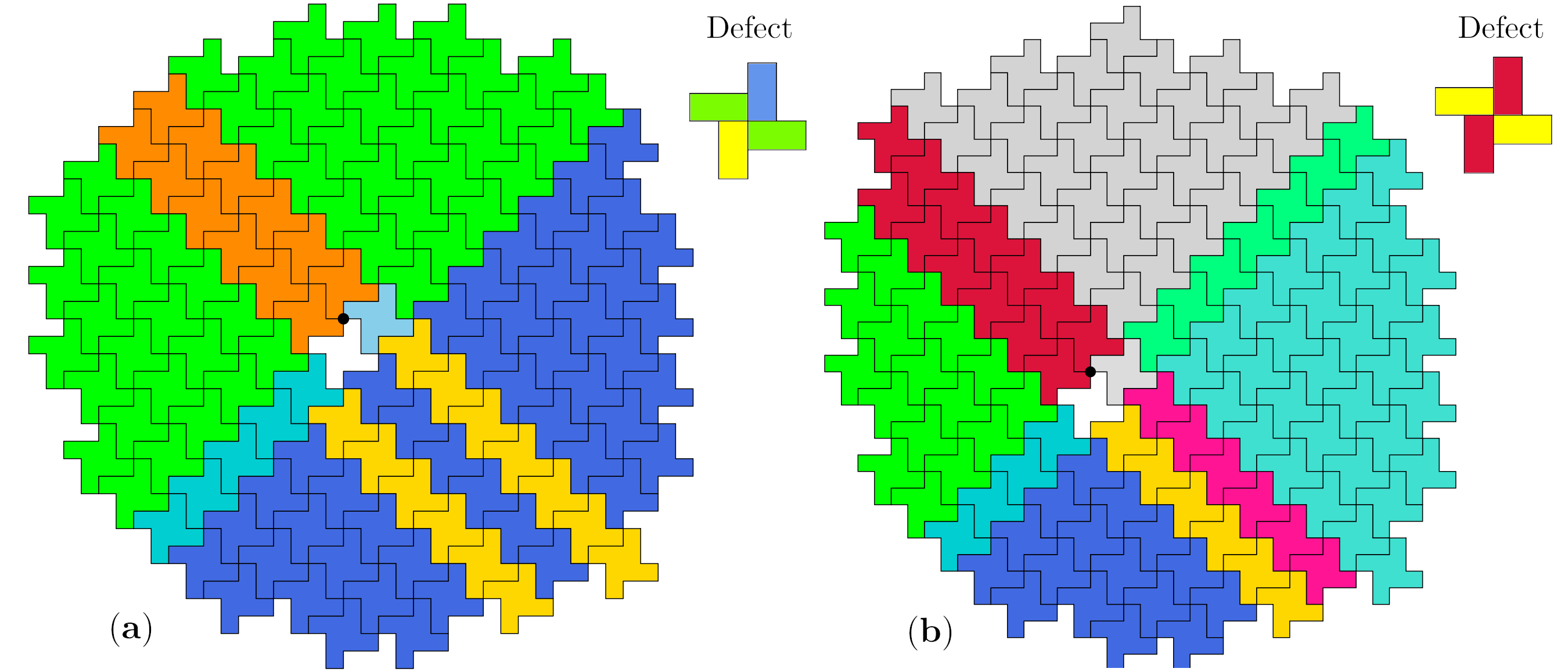} 
 	\caption{Polycrystals with  (\textbf{a}) one simple defect and (\textbf{b}) one double defect. The same 24 color ($6$ basic colors $\times4$ hues) scheme is being used. Since defects have no color assigned, they are left blank.} 
 	\label{one-defect}
 \end{figure*}

Given a finite patch $P$. Choose $M$ large enough so that all tiles from $P$ belong to $\mathcal{S}_{M}=\mathcal{S}$.  We run through all possible colorings of $\mathcal{S}$ which extend $P$ and which adhere to the adjacency rule. If the statement of the Theorem \ref{thm2} holds for any patch which consists of all tiles in $\mathcal{S}$, then the statement holds for $P$ (finite union of finite sets is finite; union of several finite and several sets of cardinality continuum has the latter cardinality).\\
\indent So, suppose we have a patch $P$ which consisted of all tiles in $\mathcal{S}$ colored. Consider all tiles which have a forced coloring. In particular, adjoin to $P$ those among them which have a common point with the square $\vert x\vert+\vert y\vert\leq M$ (like tiles \textbf{v} or \textbf{u} in Fig. \ref{saib}, if their coloring is forced by $P$). Label this extended patch as $\widetilde{P}$.\\
\indent Assume now that $P$ is not a system of $2M$ diagonal fences of length $2M$, and that the configuration $\mathcal{W}$ (see Fig. \ref{class}\textbf{B} once again) is not to be found in $\widetilde{P}$. Let at a tile $\mathbf{x}$, not belonging to the boundary of grey tiles in Fig. \ref{saib}, two diagonal fences meet (such tile necessarily exists). In other words, tiles \textbf{w} and \textbf{t} are of different colors. Then one can check that in order to avoid $\mathcal{W}$, \textbf{x} must belong to an antidiagonal fence, and also two fences around it must be built (no other variants for the coloring). Further, a light-brown staircase is then an antidiagonal fence of length $2M-1$. Moreover, tile \textbf{y} is not \textbf{R}, otherwise tiles \textbf{u} and \textbf{v} are forced and we witness $\mathcal{W}$ as a subset of $\widetilde{P}$. Consequently, a light-brown fence of length $2M-1$ extends to a fence of length $2M$. The process continues. Thus we obtain that $P$ is the union of $2M$ antidiagonal fences, each of length $2M$. \\
\indent To state what was just proved differently: if $P$ is not a part of a system of diagonal or antidiagonal type $1$ walls, $\mathcal{W}$ is always to be found in $\widetilde{P}$. We then proceed as with Theorem \ref{thm1}. This proves the first statement of the Theorem \ref{thm2}.\\
\indent It is not clear which integers can occur as orders of decent seeds. The second stated result gives a slight clue.\\
\indent Consider two red tiles given in Fig. \ref{almost}\textbf{(a)}. They are located in such a mutual position which (as the reader can check) eliminates the possibilities of type 1 or type 2 walls. From what was just proved we already know that the number of extensions is finite. Let us count.\\
\indent Consider windmills at $(1,1)$, afterwards at $(1,2)$, $(2,2)$. We see that tiles 1 and (afterwards) 2, 3 are forced to be \textbf{R}.\footnote{Just two red tiles force a countable number of \textbf{R}eds: North tile of A-mills at $(n-1,n)$, and East of C-mills at $(n,n)$, $n\in\mathbb{N}$.} Next, running the program through all admissible variants of a mill at $(0,0)$\footnote{We need to check only two cases, which can be done by hand. Indeed, permutation of colors \textbf{Y}, \textbf{B}, \textbf{G} does not impact on where program halts with a FAULT.} we get that 4 is also forced to be \textbf{R}.  Consequently, so do 5, 6. Further, tile 8 is not \textbf{R}. If a tile 7 were not \textbf{R}, then 9 and 10 would have the same (not \textbf{R}) color, and so would 11 and 12 -- a contradiction. Thus, 7 is also \textbf{R}.\\    
\indent Summing up, $2$ \textbf{R} tiles in Fig. \ref{almost}(\textbf{a}) force $7$ more \textbf{R}s shown as Fig. \ref{almost}(\textbf{b}). At this point, if an extensions exist, three remaining colors can be arbitrarily  permuted. Suppose $(13,14)=(\textbf{G},\textbf{Y})$ (Fig. \ref{almost}(\textbf{b}), the numbering of tiles is now reset). Running the program through all possible choices of tiles 1, 2, 3, 4, we obtain $6$ collections leading to perfect seeds (Fig. \ref{almost2}; the type of full field is also shown). The rest choices give a FAULT. Finally, we can fix several out of 6 variable (non-\textbf{R}) tiles in Fig. \ref{almost2} to obtain decent seeds of all orders listed. Fig. \ref{raud}(\textbf{a}) shows all \textbf{R} tiles which are forced by a seed in Fig. \ref{almost}(\textbf{a}).\\
\indent As for the seeds of arbitrarily large order, consider, for example, the one of just two \textbf{R}s shown in Fig. \ref{raud}(\textbf{b}). The following sketchy proof shows that it is decent of order at least $5\times 36$. Diagonal and antidiagonal type 1 wall system is easily ruled out. Next, suppose two tiles belong to a diagonal type 2 wall system, and brown-dotted tiles of type 2 walls lie on, say, a brown line shown in Fig. \ref{raud}(\textbf{b}). In this case dark grey tiles belong to a type 2 wall, while light-grey ones belong to a trail. A contradiction (compare to Fig. \ref{1d}(\textbf{b})). Antidiagonal type 2 wall system is also ruled out. Finally, the pattern of \textbf{R}s in Fig. \ref{raud}(\textbf{a}) can be superimposed upon our seed in 5 ways. For example, two darker-red tiles in (\textbf{a}) coincide with tiles in (\textbf{b}) if we choose to place the black dot on the left picture at the place of the larger black dot on the right. The Theorem is proved.\\   
 \indent Here is an intriguing direction of research where one can go. It is inspired by Figure 7 in \cite{socolar}, a decapod defect in the Penrose tiling which forces the full coloring. Concering a D-Cairo board, let us say that its coloring has $T$ \emph{defects}, if an adjacency rule is breached at exactly $T$ vertices. \footnote{If this rule is breached at an edge, it automatically gives $2$ ``bad" vertices.} Let us confine here to the simplest case: all possible defects occur at a vertex $(0,-1)$. And so, one (\emph{simple defect}) or both (\emph{double defect}) pairs of opposite tiles of this A-mill in the forced tiling are allowed to have the same color. Figure \ref{one-defect} shows examples for both scenarios. The supplementary file \cite{maple} contains the corresponding seeds. It would be desirable to prove an analogue of Theorem \ref{thm1} in this setting.\\
\indent One can further wonder what happens if two defects are allowed, both at odd vertices? In this case one vertex may be taken to be $(0,-1)$, while the other one can be described by an integer point $(k,\ell)$, $k+\ell$ odd. Combinatorics might get considerably more complicated and interesting.

\section{An euclidean model}
\label{toy}
\begin{figure}
	\begin{center}
	\includegraphics[scale=0.33]{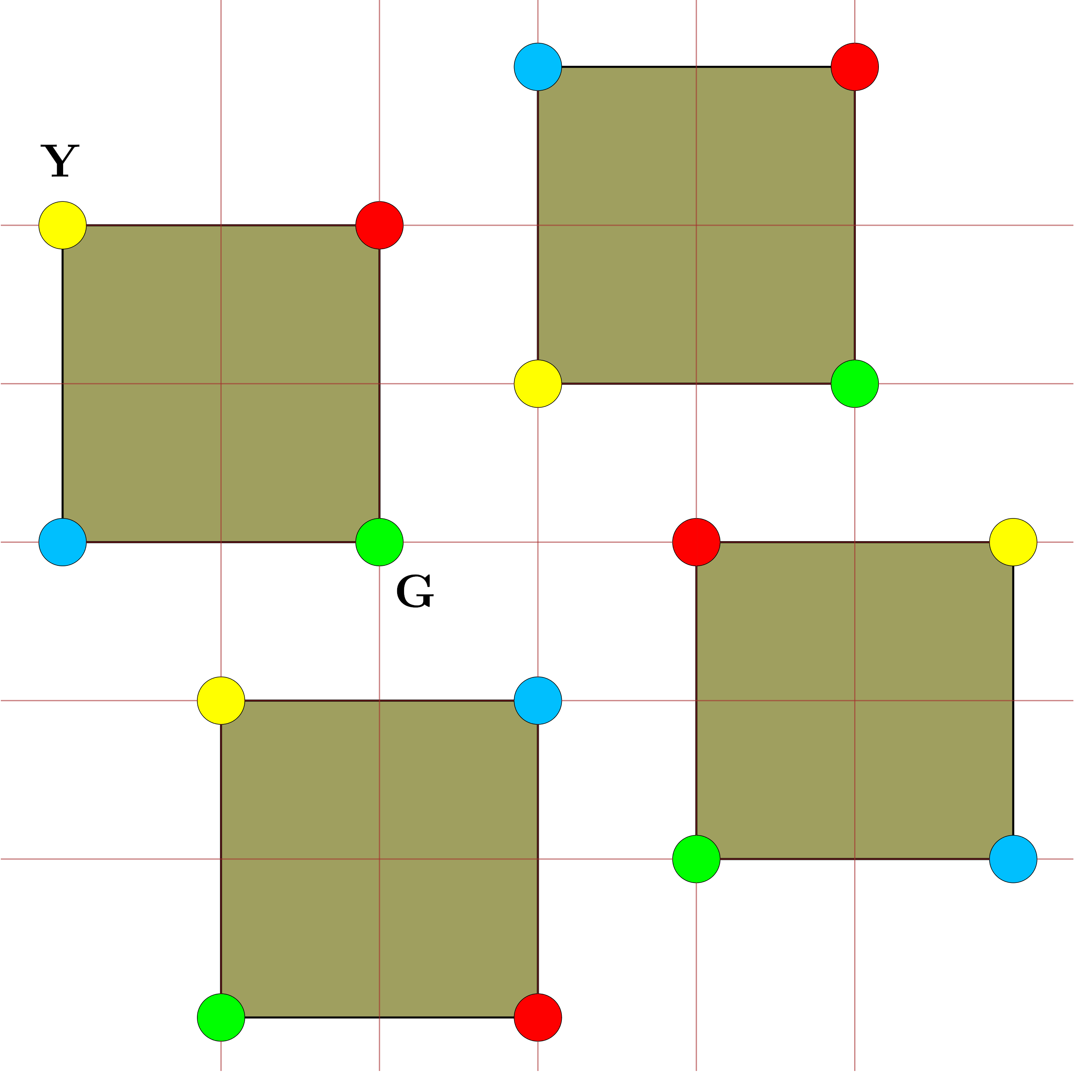} 
	\caption{$4$ molecules corresponding to $4$ central windmills in Figure \ref{trys} (top left)}.
	\label{crys} 
	\end{center}
\end{figure}

\indent Let us consider our tiles as $24$ versions of squares (see Fig. \ref{unit} again). We will show that the matching rules can be realised by euclido-geometric restrictions.  Fig. \ref{crys} shows a section of a model, which we will now explain.\\   
\indent  For this purpose consider a $2\times 2\times 2$ cube (``a molecule"). It consists of $8$ ``atoms", two copies of each of \textbf{R}, \textbf{G}, \textbf{B}, \textbf{Y}, identical atoms being placed on opposite ends of long diagonals. All molecules thus only differ in their spacial orientation. Inside the $xy$ plane, consider the lattice $\mathbb{Z}^{2}$ and its index $10$ sublattice $\mathcal{L}=(-1,3)\mathbb{Z}\times(3,1)\mathbb{Z}$. For each $\mathbf{p}\in\mathcal{L}$, let us place a molecule with $\mathbf{p}$ as its centre, so that the upper face of a molecule is parallel to the $xy$ plane, and its edges are aligned with the axes $x,y$.  Figure \ref{crys} shows $4$ specimens.  The process now runs as follows: at every discrete tick of a time, one molecule makes several $\pm 90^{\circ}$ rotations around its centre with respect to any of the axes $x$, $y$, and $z$ (in the picture, the latter is directed towards the viewer). This is called \emph{an orientation change}.
The adjacency rule translates into the language of molecules as follows: 
 \begin{center}
 	\underline{Rule}:\emph{ the distance between atoms of the same kind\\
 		must be greater than} $2$.
 \end{center}
Note that for each atom there are exactly $7$ other atoms at a distance $\leq 2$. For example, if we choose \textbf{G} as shown in Figure \ref{crys}, $6$ of these closest neighbours belong to the plane $z=1$, while the remaining one, namely, \textbf{Y}, lies on the plane $z=-1$ immediately below \textbf{G}. The former is paired with \textbf{Y} shown on the picture. During the process, each molecule strives to attain orientation which is compliant with the \underline{Rule}. According to R. Penrose, Euclidean geometry is the first SUPERB physical theory \cite{penrose}. Thus, the above construction gives a neat ``physical model" of our coloring problem, translating matching rules into a single geometric restriction. Possibly, this toy is implementable relying solely on mechanics (like an analougue of Rubik's cube for each molecule) and magnetism (a pair of magnets for each atom). Six faces of a cube can be given 6 colours and certain textures, which, if rotated, would reflect a slightly different hue. This would be in correspondence with coloring scheme used in two previous sections.   

\section{Type $7$ pentagonal tiling}
\begin{figure}
	\begin{center}
	\includegraphics[scale=0.45]{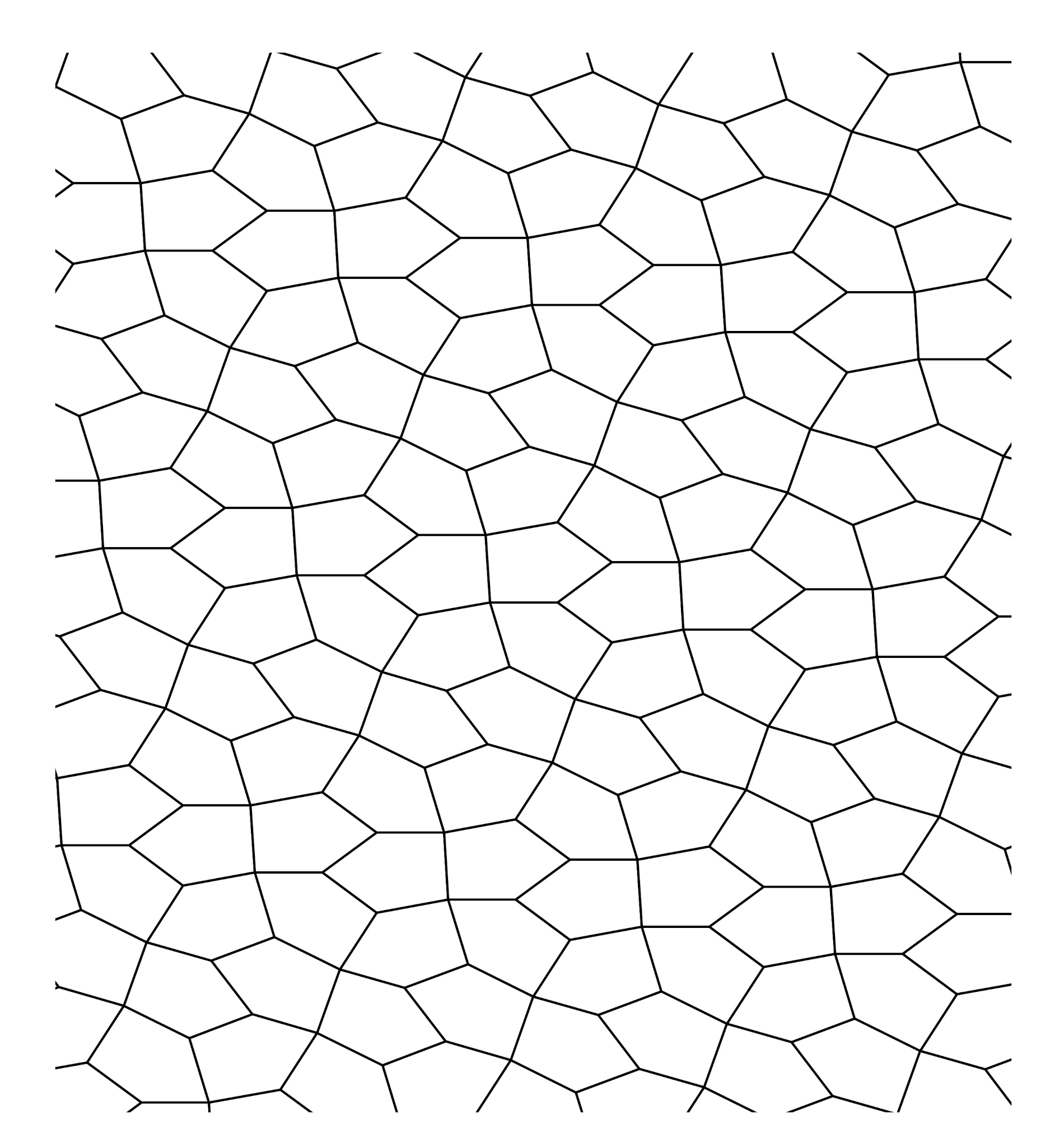} 
	\caption{Pentagonal type $7$ tiling} 
	\label{sept}
	\end{center}
\end{figure}

\begin{figure}
	\begin{center}
	\includegraphics[scale=0.30]{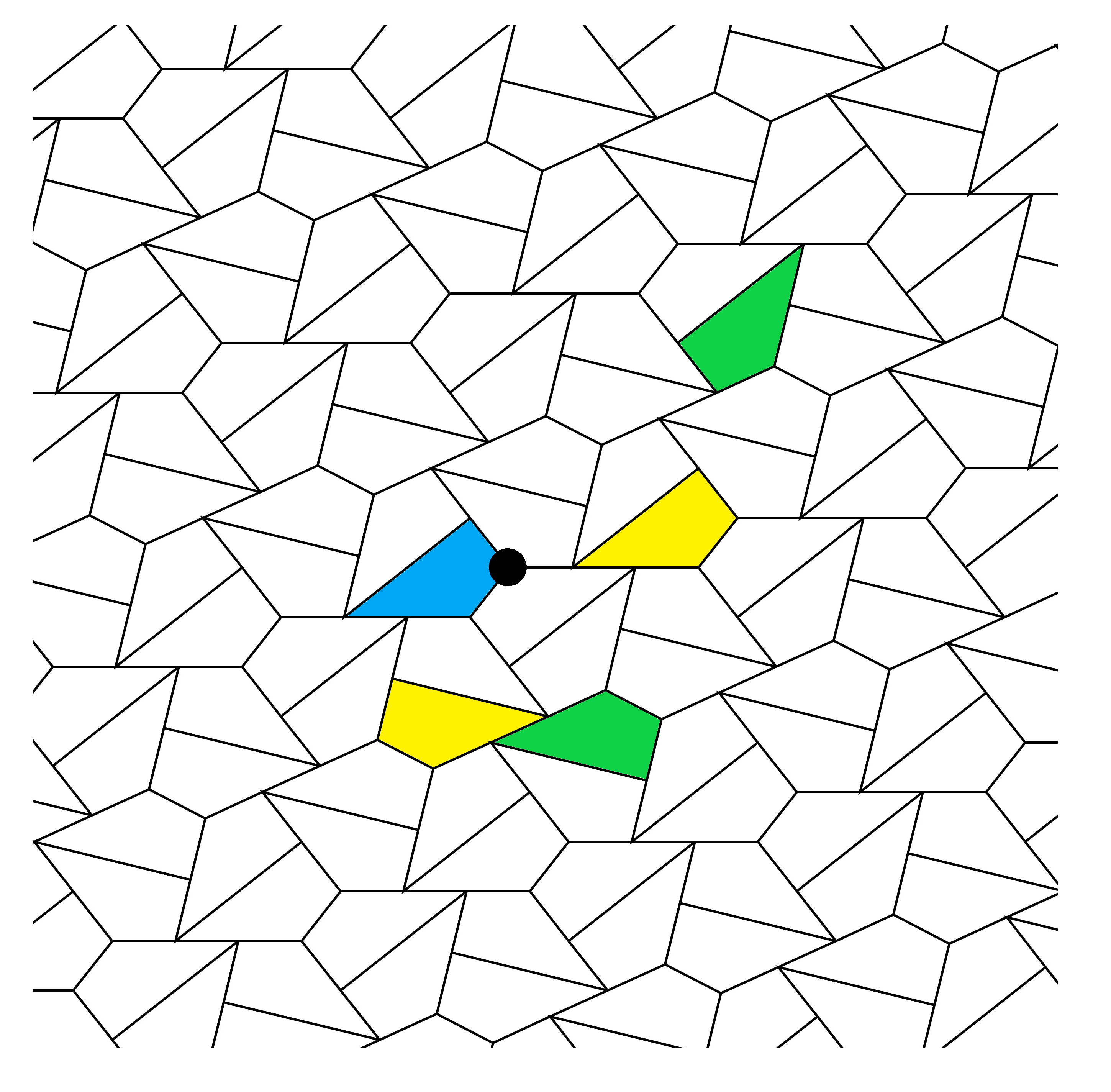} 
	\caption{A perfect seed} 
	\label{sept2}
	\end{center}
\end{figure}

\begin{figure}
	\begin{center}
		\includegraphics[scale=0.40]{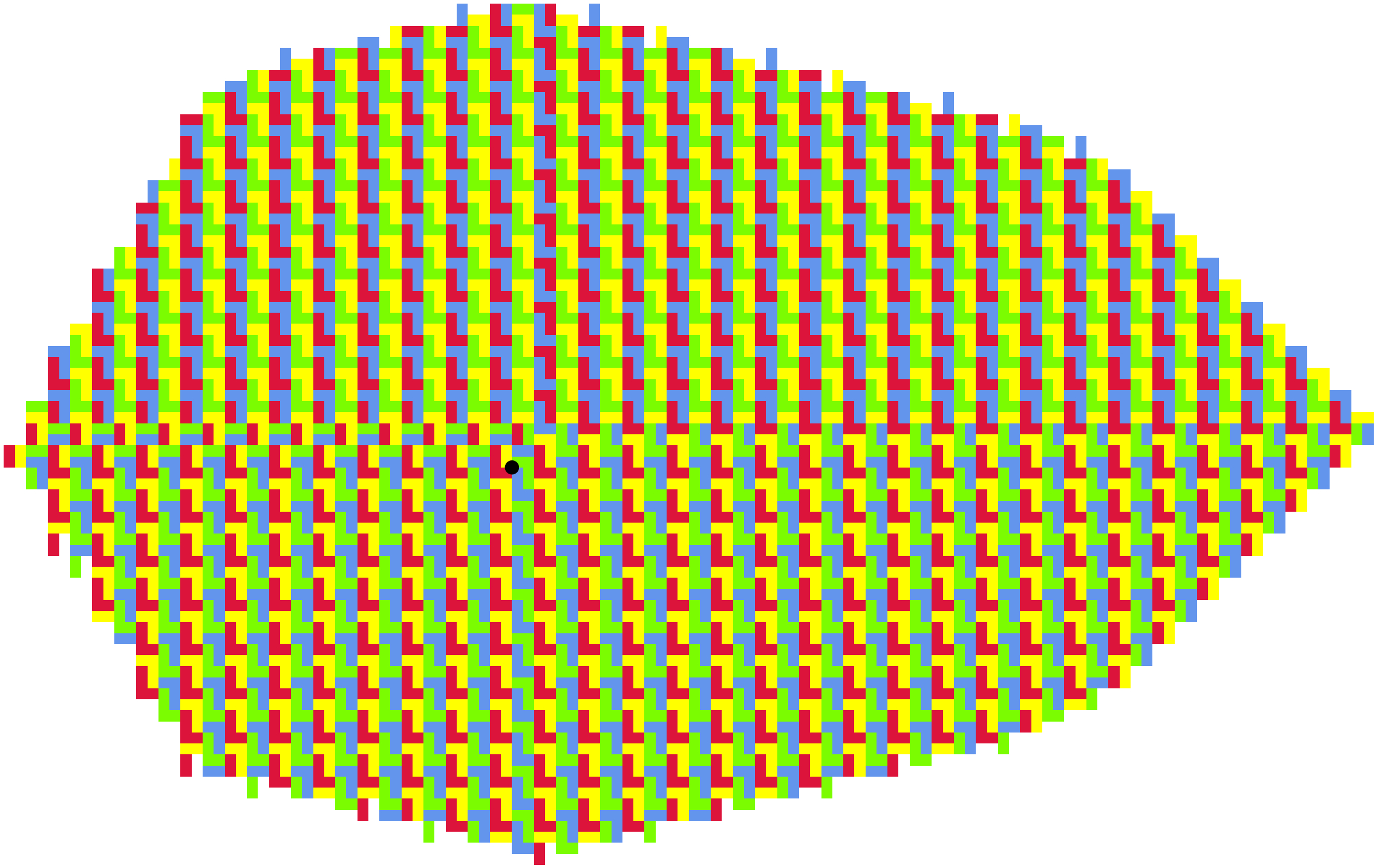} 
		\caption{$80$ generations of the seed in Fig. \ref{sept2}} 
		\label{80-penta}
	\end{center}
\end{figure}

\begin{figure}
	\includegraphics[scale=0.42]{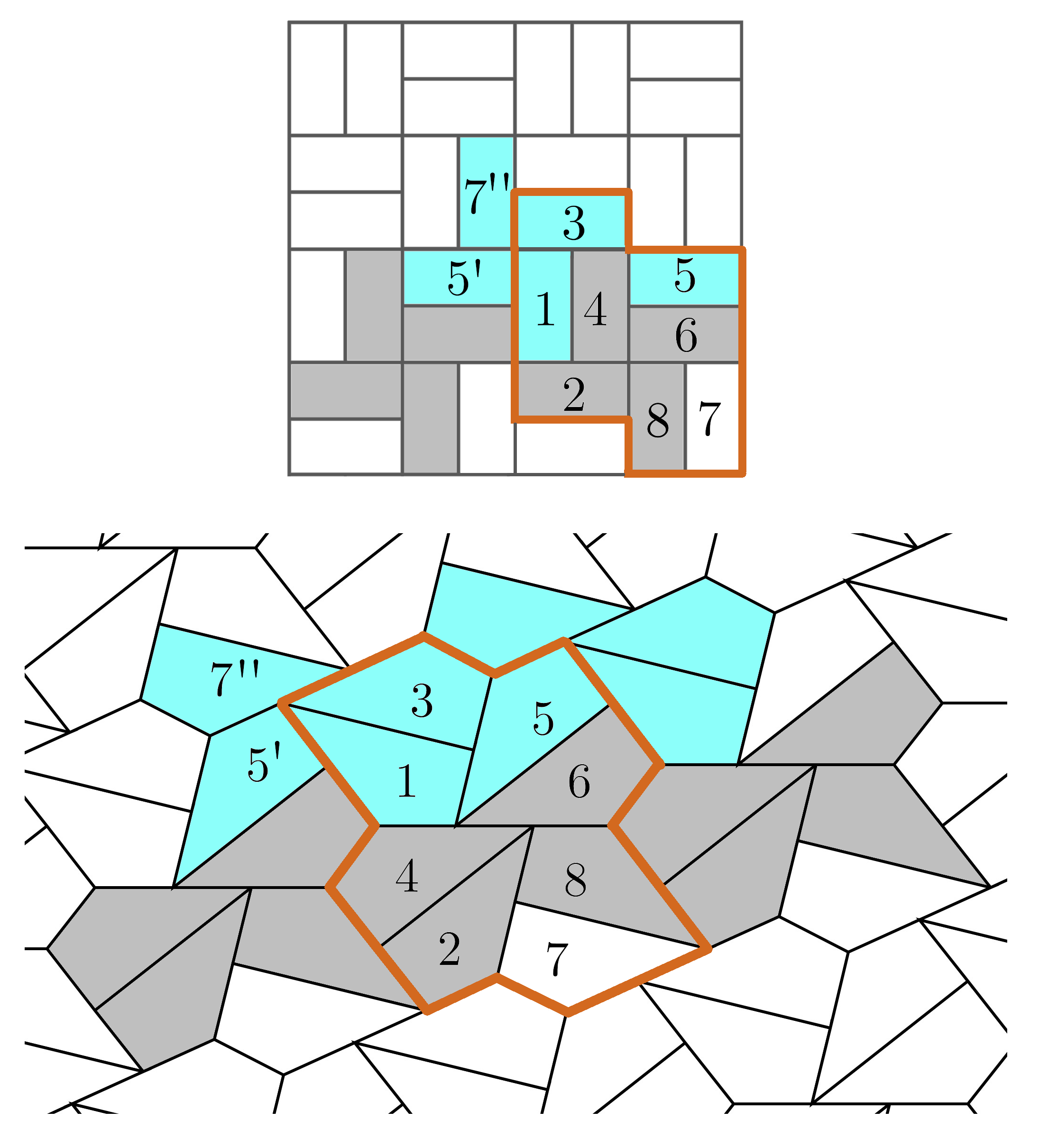} 
	\caption{The map from type 7 pentagonal to the D-Cairo} 
	\label{septyni}
\end{figure}

\begin{figure*}
	\begin{center}
	\includegraphics[scale=0.70]{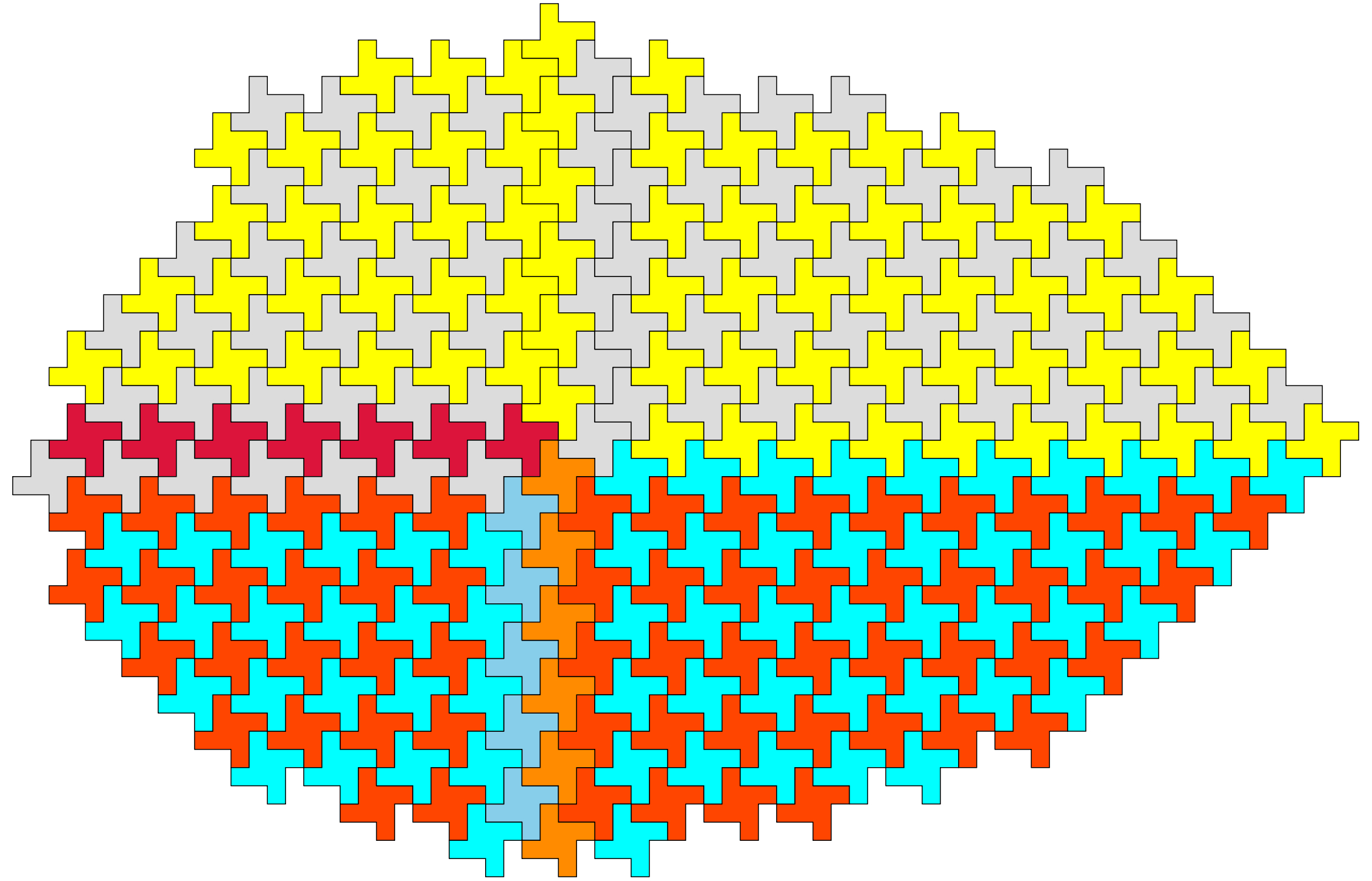} 
	\caption{``Even" and ``odd" part of the same field} 
	\label{septyni-f}
	\end{center}
\end{figure*}

To show that the D-Cairo tiling is not an isolated phenomenon where forced grain boundary occurs, consider $2$-isohedral pentagonal tiling of type $7$ (discovered by R. Kershner, Figure \ref{sept}). Its extreme case is topologically equivalent quadrilateral tiling (denoted in this section by K, see Figure \ref{sept2}).\footnote{One can play around all $15$ known types of pentagonal tilings on a Wolfram tool \cite{wolfram}.} \\
\indent Perfect seeds do exist for the tiling K, too. One of them is shown in Figure \ref{sept2}. The fundamental domain (shown in Figure \ref{septyni} as a brown contour) consists of 8 tiles. There appears to be $196$ types of colorings of 8 tiles, which comply to the adjacency rule. One can invent a new coloring scheme with $196$ colors in the palette and present the pictures obtained.\\
\indent But none of that is needed. In fact, it is possible to work with the same D-Cairo tiling. To demonstrate this, let us map all tiles 1through 8 in the fundamental domain of a K-tiling to the properly chosen 8 tiles in D-Cairo tiling, as shown in Fig. \ref{septyni}. Thus, the same D-Cairo geometry is being used, only the set of $7$ neighbours for each tile of the field is redefined.  \\
\indent Now, if we calculate $80$ generations of the perfect seed in Fig. \ref{sept2} and transport the picture to the D-Cairo field via a map given, we obtain Fig. \ref{80-penta}.\\
\indent One clearly sees a grain boundary emerging. What is unexpected is that the growth of the polycrystal in the direction to the readers right hand is faster that to the left. Let us take a final look to the map in Fig. \ref{septyni}. Tiles 2, 4, 6, 8 from one fundamental domain mutually touch one another. We can color the corresponding C-mill in the D-Cairo field according to our previous coloring convention. The same applies to tiles 1, 3, 5', 7'', this time tiles 5' and 7'' come from two neighbouring fundamental domains. The picture obtained is shown in Fig. \ref{septyni-f}.\\ 
\indent In order to finally convince oneself that this grain boundary is of a different kind than those described in Section \ref{all-poly}, let us start again from the seed in Fig. \ref{sept2} and count the number (say, $k(n)$) of tiles we are able to color at a step $n$. This produces the sequence
\begin{eqnarray*}
	1, 3, 3, 3, 2, 2, 4, 3, 4, 6, 8, 9, 11, 12, 14, 12, 16, 16,\ldots
	\end{eqnarray*}  
It appears that $k(n+1)-k(n)$ is also eventually periodic. Here is the final result: for $n\geq 37$, $k(n)=\big{\lfloor}\frac{15n}{14}\big{\rfloor}+c(n)$, where $c(n)$ is a sequence with period $336$. Proving exact values for coordination sequences even for the basic tilings is not a simple task (see \cite{sloan}). Thus, to rigorously prove the results stated about the sequences $t(n)$ and $k(n)$ might be much harder than to check this experimentally.  \\
\indent Going over many other examples given, say, in the masterpiece of the field \cite{shephard}, will hopefully produce even more grain boundaries. It is plausible that the technique we used is applicable in a much general setting. That is, first one needs to understand what are periodic walls. Second, find a starting seed(s) which always occur(s) provided periodic walls are absent. Third, find a well-chosen set of tiles and perform the full computer-assisted check, discovering other types of walls on the way. \\
\indent\emph{Acknowledgement.} The author sincerely thanks the anonymous referee for benevolence and detailed remarks.   



\end{document}